\documentclass[aap,MSNbibl,seceqn,dvips]{arximspdf}
\usepackage{mathbh}
\usepackage{graphicx}
\doi{10.1214/12-AAP862} 
\volume{23}
\issue{3}
\pubyear{2013}
\firstpage{987}
\lastpage{1024}

\makeatletter
\newcommand{\pmb}[1]{\bolds{#1}}
\newcommand{\pmbb}[1]{\mathbf{#1}}
\newcommand{\eqref}[1]{(\ref{#1})}
\newtheorem{theorem}{Theorem}[section]
\newtheorem{lemma}[theorem]{Lemma}
\newproclaim{definition}[theorem]{Definition}
\newtheorem{corollary}[theorem]{Corollary}
\newproclaim{remark}{Remark}

\newcommand{\R}{\mathbb{R}}
\newcommand{\C}{\mathbb{C}}

\newcommand{\Z}{\mathbb{Z}}
\newcommand{\N}{\mathbb{N}}

\newcommand{\var}{\operatorname{Var}}

\newcommand{\PT}{\mathbb{P}_{\Theta}}

\newcommand{\la}{\lambda}
\newcommand{\tr}{\operatorname{Tr}}
\newcommand{\one}{\mathbh{1}}
\newcommand{\SN}{\mathfrak{S}_N}

\renewcommand{\Re}{\operatorname{Re}}

\makeatother

\begin{document}
\begin{frontmatter}

\title{Random permutation matrices under the generalized Ewens measure}
\runtitle{Random permutation matrices}

\begin{aug}
\author[A]{\fnms{Christopher}~\snm{Hughes}\ead[label=e1]{christopher.hughes@york.ac.uk}},
\author[B]{\fnms{Joseph}~\snm{Najnudel}\ead[label=e2]{joseph.najnudel@math.uzh.ch}},
\author[B]{\fnms{Ashkan}~\snm{Nikeghbali}\ead[label=e3]{ashkan.nikeghbali@math.uzh.ch}}
\and
\author[C]{\fnms{Dirk} \snm{Zeindler}\corref{}\thanksref{T1}\ead[label=e4]{zeindler@math.uni-bielefeld.de}}
\thankstext{T1}{Supported by the Swiss National Science Foundation (SNF).}
\runauthor{Hughes, Najnudel, Nikeghbali and Zeindler}
\affiliation{University of York, Universit\"at Z\"urich, Universit\"at
Z\"urich and~Universit\"at~Bielefeld}
\address[A]{C. Hughes\\
Department of Mathematics\\
University of York\\
York YO10 5DD\\
United Kingdom\\
\printead{e1}}

\address[B]{J. Najnudel\\
A. Nikeghbali\\
Institut f\"ur Mathematik\\
Universit\"at Z\"urich\\
Winterthurerstrasse 190\\
8057-Z\"urich\\
Switzerland\\
\printead{e2}\\
\phantom{E-mail: }\printead*{e3}}

\address[C]{D. Zeindler\\
Fakult\"at f\"ur Mathematik\\
Universit\"at Bielefeld\\
Sonderforschungsbereich 701\\
Postfach 10 01 31\\
33501 Bielefeld\\
Germany\\
\printead{e4}}

\end{aug}

\received{\smonth{11} \syear{2011}}
\revised{\smonth{3} \syear{2012}}

%
\begin{abstract}
We consider a generalization of the Ewens measure for the
symmetric group, calculating moments of the characteristic
polynomial and similar multiplicative statistics. In addition,
we study the asymptotic behavior of linear statistics (such as
the trace of a permutation matrix or of a wreath product) under
this new measure.
\end{abstract}

%


\end{frontmatter}

\section{Introduction}\label{sec1}
The trace of a matrix is one of the most natural additive class
functions associated to the spectra of a matrix. Traces of
unitary matrices chosen randomly with Haar measure have been
much studied, for example by Diaconis and Shahshahani
\cite{DiacShah94}, Diaconis and Evans~\cite{DiaEva01} and Rains
\cite{Rains97} using methods from representation theory, by
Diaconis and Gamburd~\cite{DiaGam04} using combinatorics and
using methods from mathematical physics by Haake et
al.~\cite{Haa96}.

Another natural class function, this time multiplicative, is
the characteristic polynomial, and the distribution of
characteristic polynomials of random unitary matrices has been
studied by many authors, including Keating and Snaith~\cite{KS00a}
and Hughes, Keating and O'Connell~\cite{HKO01}.

From the characteristic polynomial one can find the number of
eigenvalues lying in a certain arc (since the underlying matrix
is unitary, all the eigenvalues lie on the unit circle). The
problem of studying the number of eigenvalues lying in an arc
was studied by Rains~\cite{Rains97} and Wieand~\cite{Wie02}
who found a very interesting correlation structure when
multiple arcs were considered, and Hughes, Keating and
O'Connell~\cite{HKO01} who made the connection with
characteristic polynomials.

One of the reasons for such an extensive study into random
unitary matrices and their spectra is that the statistical
distribution of the eigenvalues is expected to have the same
behavior as the zeros of the Riemann zeta function (see
Montgomery~\cite{Mon73} and Keating and Snaith~\cite{KS00a}).

The statistics of the distribution of spectra of infinite
subgroups of the unitary group, such as the symplectic and
orthogonal groups, are also expected to model families of other
$L$-functions (Keating and Snaith~\cite{KS00b}), and the
distribution of traces for these groups have been studied,
frequently in the same papers.

However, the statistics of the spectra of finite subgroups of
the unitary group, such as the permutation group, is not so
well studied, though there are many results known.

Wieand~\cite{Wie00} studied the number of eigenvalues of a
uniformly random permutation matrix lying in a fixed arc, and
Hambly et al.~\cite{HKOS} found
corresponding results for the characteristic polynomial, making
the same connection between the characteristic polynomial and
the counting function of eigenvalues.

In all these cases, the permutation matrices were chosen with
uniform measure and the results were similar to those found for
the full unitary group with Haar measure. However, there were
some significant differences primarily stemming from the fact
that the full unitary group is rotation invariant, so the
characteristic polynomial is isotropic. The group of
permutation matrix is clearly not rotation invariant, and the
distribution of the characteristic polynomial depends weakly on
the angle of the parameter. Most results require the angle to
have the form $2\pi\tau$ with $\tau$ irrational and of finite
type. (It is worth pointing out that those angles which are not
of this type, have Hausdorff dimension zero). More recently Ben
Arous and Dang~\cite{benarousdang} have extended some of the
results of Wieand to more general measures together with new
observations very specific to permutation matrices. In
particular, they prove that the fluctuations of smooth linear
statistics (with bounded variance) of random permutation
matrices sampled under the Ewens measure are asymptotically
non-Gaussian but infinitely divisible.

A permutation matrix has the advantage that the eigenvalues are
determined by the cycle-type $\la= (\la_1,\ldots,\la_\ell)$ of
the corresponding permutation. This allows us to write down in
many situations explicit expression for the studied object. For
example, the characteristic polynomial of an $N\times N$
permutation matrix $M$ is
%
\begin{eqnarray} \label{eqdefZnintro}
Z_N(x) &=& \det(I-x M) = \prod_{j=1}^N (1-x
e^{i\alpha_j}) = \prod_{m=1}^{\ell(\lambda)}
(1-x^{\lambda_m})
\nonumber
\\[-8pt]
\\[-8pt]
\nonumber
 &=& \prod_{k=1}^N
(1-x^k )^{C_k},
\end{eqnarray}
where $e^{i\alpha_1}, \ldots, e^{i\alpha_N}$ are the
eigenvalues of $M$ and $C_k$ is the number of cycles of length
$k$ (i.e., the number of $j$ such that $\lambda_j=k$).

Equation \eqref{eqdefZnintro} has been used by Dehaye and
Zeindler~\cite{Zei10} to introduce multiplicative class
functions (i.e., invariant under conjugation) associated to function~$f$. These functions have the
form
%
\begin{equation}
W^N(f)(x) = \prod_{k=1}^N f(x^k)^{C_k}
\label{eqdefWnintro}
\end{equation}
and generalize the characteristic polynomial, which is the case
$f(x) = (1-x)$. This nice generalization is not possible for
the unitary group because the eigenvalues lack the structure of
eigenvalues of permutation matrices that allows
\eqref{eqdefZnintro} to hold. The most natural analogue for
unitary matrices is Heine's identity which connects Haar
averages of multiplicative class functions with determinants of
Toeplitz matrices.

Equation \eqref{eqdefZnintro} can also be used to express linear
statistics or traces of functions of random permutation matrices in
terms of the cycle counts. More precisely, if we identify a random
permutation matrix $M$ with the permutation $\sigma$ it represents, we
have the following definition.
%
\begin{definition}
\label{deftraceF}
Let $F\dvtx S^1 \to\C$ be given. We then
define the trace of $F$ to be the function $\tr(F)\dvtx\SN\to\C$
with
%
\begin{equation}\label{eqbahroumi}
\tr(F)(\sigma):= \sum_{k=1}^N F(\omega_k),
\end{equation}
where $(\omega_k)_{k=1}^N$ are the eigenvalues of $\sigma$ with
multiplicity.
\end{definition}

Observe that when $F(x) = x^d$, we have $\tr(F)(\sigma) =
\tr(\sigma^d)$, and this justifies the use of the terminology
trace. The trace of a function is also referred to as a linear
statistic on $\SN$.

%
\begin{lemma}
\label{lemtracefunctionexplicit}
Let $F\dvtx S^1 \to\C$ and $\sigma\in\SN$ with cycle type $\la$
be given. We then have
%
\begin{equation}
\tr(F)(\sigma)
=
\sum_{k=1}^{N} kC_k \Delta_{k}(F)
\end{equation}
with $\Delta_{k}(F):= \frac{1}{k} \sum_{m=1}^k F( e^{2\pi
i m/k} )$.
\end{lemma}

\begin{pf} This follows immediately from equation
\eqref{eqdefZnintro}.
\end{pf}

The reason why the expressions in \eqref{eqdefZnintro} and
\eqref{eqdefWnintro} are useful is that many things are
known about the cycle counts $C_k$. For example, the cycle
counts $C_k$ converge weakly to independent Poisson random
variables $P_k$, with mean $1/k$. Also useful in this context
is the Feller coupling, since in
many situations this allows one to
replace~$C_k$ by $P_k$. Several details on the cycle counts and
the Feller coupling can be found in the book by
Arratia, Barbour and Tavar{\'e}~\cite{ABT02}.

These results all concern the uniform measure, where each
permutation has weight $1/N!$, or the
Ewens measure, where the probability is proportional to the
total number of cycles, $\theta^{\ell(\lambda)} / (N!h_N)$,
with $h_N$ the required normalization constant (the case $\theta=1$
corresponding to the uniform measure). A common ingredient in all the
above cited works is the use of the Feller coupling (details can again
be found in the book by Arratia, Barbour and
Tavar{\'e}~\cite{ABT02}) and some improvements on the known bounds for
the approximation given by this coupling (see~\cite{benarousdang},
Section~4).

In recent years, there have been many works in random matrix theory
aimed at understanding how much the spectral properties of random
matrices depend on the probability distributions of its entries. Here a
similar question translates into how are the linear and multiplicative
(i.e., multiplicative class functions) statistics affected if one
considers more general probability distributions than the Ewens measure
on the symmetric group? The Ewens measure can be naturally generalized
to a weighted
probability measure which assigns to the permutation matrix $M$ (i.e.,
to the associated permutation) the
weight
\[
\frac{1}{N! h_N} \prod_{k=1}^N \theta_k^{C_k},
\]
where $h_N$ is a normalization constant. The Ewens measure
corresponds to the special case where $\theta_k = \theta$ is a
constant. This measure has recently appeared in mathematical
physics models (see, e.g.,~\cite{BUV} and
\cite{cyclestructureueltschi}) and one has only recently
started to gain insight into the cycle structures of such
random permutations. One major obstacle with such measures is
that there exists nothing such as the Feller coupling and
therefore the classical probabilistic arguments do not apply
here. In a recent work, Nikeghbali and Zeindler~\cite{NikZei}
propose a new approach based on combinatorial arguments and
singularity analysis of generating functions to obtain
asymptotic expansions for the characteristic functions of the
$C_k$'s as well as the total number of cycles, thus extending
the classical limit theorems (and some distributional
approximations) for the cycle structures of random permutations
under the Ewens measure. In this paper we shall use the methods
introduced in~\cite{NikZei}, namely, some combinatorial lemmas,
generating series and singularity analysis to study linear
and multiplicative statistics for random permutation matrices
under the general weighted probability measure. In fact, we
shall consider the more general random matrix model obtained
from the wreath product $S^1 \wr\SN$ (see, e.g.,~\cite{Wie03}); this amounts to replacing the $1$'s in the
permutation matrices by independent random variables taking
values in the unit circle $S^1$.
The distribution of eigenvalues of such
matrices (alongside other generalizations) has been studied previously by Najnudel and
Nikeghbali~\cite{NajNik11}.
It should be noted that many
groups closely related to $\SN$ exhibit such matrices, for
instance, the Weyl group of $\operatorname{SO}(2N)$.

More precisely this paper is organized as follows.

In Section~\ref{secpreparation} we fix some notation and
terminology, recall some useful combinatorial lemmas together
with some results of Hwang (and some slight extensions) on
singularity analysis of generating functions. In particular, we
shall introduce two relevant classes of generating functions
according to their behavior near singularities on the circle
of convergence. In this article, we shall state our theorems
for random matrices under the generalized Ewens measures for
which the generating series of $(\theta_k)_{k\geq1}$ is in one
of these two classes.

In Section~\ref{sectACF} we study the multiplicative class
functions associated to a function~$f$ and obtain the
asymptotic behavior of the joint moments. In particular, we
extend earlier results of~\cite{Zei10,Zei10a} and of
\cite{HKOS} on the characteristic polynomial of uniformly
chosen random permutation matrices.

In Section~\ref{sectraces} we focus both on the traces of
powers and powers of traces which are classical statistics in
random matrix theory. In fact we prove more generally that the
fluctuations of the linear statistics for Laurent polynomials
are asymptotically infinitely divisible (they converge in law
to an infinite weighted sum of independent Poisson variables).
We also establish the convergence of the integer moments of
linear statistics for functions of bounded variation together
with the rate of convergence.

In Section~\ref{secgeneralcase} we consider the more general
model consisting of the wreath product $S^1 \wr\SN$ and study
the linear statistics for general functions $F$ in
\eqref{eqbahroumi}. In such models, the $1$'s in the
permutation matrix are replaced with $(z_j)_{1\leq j\leq N}$
which are i.i.d. random variables taking their values on the
unit circle $S^1$. In this framework, Lemma
\ref{lemtracefunctionexplicit} can be naturally extended
(see Lemma~\ref{lembahroumi}) and the quantity
%
\begin{equation}
\Delta_k(F,z) := \frac{1}{k} \sum_{\omega^k = z} F(\omega)
\end{equation}
naturally appears in our technical conditions. Under some
conditions on rate of convergence to $0$ of the $L^1$-norm of $
\Delta_k(F,z)$, and some assumptions on the singularities of
the generating series of $(\theta_k)_{k\geq1}$, we are able to
compute the asymptotics of the characteristic function of
$\tr(F)$ with a good error term. From these asymptotics we are
able to compute the fluctuations of $\tr(F)$. We also translate
our conditions in terms of the Fourier coefficients of $F$,
where $F$ has to be in some Sobolev space $H^s$.

In Section~\ref{sectEwensDivVar} we still work within the
framework of the wreath product $S^1 \wr\SN$ and consider the
case where the variance of $\tr(F)$ is diverging. This time we
restrict ourselves to the Ewens measure since our methods do
not seem to apply in this situation. Hence, we go back to
probabilistic arguments (i.e., use the Feller coupling) to prove
that under some technical conditions on $F$, the fluctuations
are Gaussian. In fact, we essentially adapt the proof by Ben
Arous and Dang~\cite{benarousdang} to these more general
situations. Nonetheless our theorem is also slightly more
general in that it applies to a larger class of functions $F$.

\section{The generalized Ewens measure, generating series and
singularity analysis}
\label{secpreparation}
In this section we fix notation and we recall some facts about the
symmetric group and generating functions, as well as the main results
from singular analysis (we also provide some variants and extensions
for the purpose of this paper). Our presentation closely follows~\cite{NikZei}.
\subsection{Some combinatorial lemmas and the generalized Ewens measure}
We present in this section some basic facts about $\SN$ and
then define the generalized Ewens measure. We give here only a
very short overview and refer to~\cite{ABT02} and~\cite{macdonald} for more details.

\subsubsection{Conjugation classes and functions on $\SN$}
\label{secconjclassofSn} We first take a closer look at the
conjugation classes of the symmetric group $\SN$ (the group of
all permutations of a set of $N$ objects). We only need to
consider the conjugation classes since all probability measures
and functions considered in this paper are invariant under
conjugation (i.e.,~they are class functions). It is well known
that the conjugation classes of $\SN$ can be parameterized with
partitions of $N$.
%
\begin{definition}
\label{defpart}
A partition $\la$ is a sequence of nonnegative integers $\la_1
\ge\la_2 \ge\cdots$ eventually trailing to $0$'s, usually
omitted. The size of the partition is $|\lambda|:= \sum_m
\lambda_m$. We call $\la$ a partition of $N$ if $|\la| = N$,
and this will be denoted by $\la\vdash N$. The length of $\la$
is the largest $\ell$ such that $\la_\ell\ne0$.
\end{definition}

Let $\sigma\in\SN$ be arbitrary. We can write $\sigma=
\sigma_1\cdots\sigma_\ell$ with $\sigma_m$, $1\leq m\leq\ell$,
disjoint cycles of length $\la_m$. Since disjoint cycles
commute, we can assume that $\la_1\geq\la_2\geq\cdots
\geq\la_\ell$. We call the partition
$\la=(\la_1,\la_2,\ldots,\la_\ell)$ the \textit{cycle-type} of~$\sigma$.
We write $\mathcal{C}_\la$ for the set of all $\sigma
\in\SN$ with cycle type~$\la$. One
can now show that two elements $\sigma,\tau\in\SN$ are
conjugate if and only if $\sigma$ and $\tau$ have the same
cycle-type and that $C_\la$ are the conjugation classes of
$\SN$. Since this is well known, we omit the proof and refer the
reader to
\cite{macdonald} for more details.

\begin{definition}
\label{defCmandK0n} Let $\sigma\in\SN$ be given with
cycle-type $\la$. The cycle numbers~$C_k$ are defined as
%
%
\begin{equation}\label{cyclebahroumi}
C_k = C_k(\sigma) := \#\{m \dvtx\la_m = k\}
\end{equation}
and the total number of cycles $T(\sigma)$ is
%
\begin{equation}
T(\sigma):= \sum_{k=1}^N C_k.
\end{equation}
\end{definition}

The functions $C_k(\sigma)$ and $T(\sigma)$ depend only on the
cycle type of $\sigma$ are thus class functions. Clearly
$T(\sigma)$ equals $\ell(\lambda)$, the length of the partition
corresponding to~$\sigma$.

All expectations in this paper have the form $\frac{1}{N!}
\sum_{\sigma\in\SN} u(\sigma)$ for a class function~$u$. Since
$u$ is constant on conjugation classes, it is more natural to
sum over all conjugation classes. We thus need to know the size
of each conjugation class.

%
\begin{lemma}
\label{lemsizeofconjclasses}
We have
%
%
\begin{equation}
|\mathcal{C}_\la|=\frac{|\SN|}{z_\la}\qquad
\mbox{with }
z_\la:=\prod_{k=1}^{N} k^{C_k}C_k!
\end{equation}
with $C_k$ defined in \eqref{cyclebahroumi}, and
%
\begin{equation}\label{eqsumclassfuncs}
\frac{1}{N!} \sum_{\sigma\in\SN} u(\sigma)
=
\sum_{\la\vdash N}\frac{1}{z_\la}u(\mathcal{C}_\la)
\end{equation}
for a class function $u\dvtx\SN\to\C$.
\end{lemma}

\begin{pf}
The first part can be found in~\cite{macdonald} or in~\cite{bump},
Chapter 39.
The second part follows immediately from the first part.
\end{pf}

\subsubsection{Definition of the generalized Ewens measures}
\label{secdefofweightedprob}
We now define the generalized Ewens measures.
%
\begin{definition}
\label{defweightedprobabililtymeasure}
Let $\Theta= (\theta_k )_{k=1}^\infty$ be a
sequence of strictly positive numbers. We define for
$\sigma\in\SN$ with cycle-type $\la$
%
\begin{equation}
\PT[\sigma]
:=
\frac{1}{h_N N!} \prod_{m=1}^{\ell(\la)} \theta_{\la_m}
=
\frac{1}{h_N N!} \prod_{k= 1}^N \theta_k^{C_k(\sigma)}
\label{eqPThetawithpartition}
\end{equation}
with $h_N = h_N(\Theta)$ a normalization constant and $h_0:=1$.
\end{definition}

The second equality in \eqref{eqPThetawithpartition} follows
immediately from the definition of $C_k$
(Definition~\ref{defCmandK0n}). The uniform measure and the
Ewens measure are special cases, with $\theta_k \equiv1$ and
$\theta_k\equiv\theta$ a constant, respectively.\vadjust{\goodbreak}

We now introduce two generating functions closely related to
$\PT$:
%
\begin{equation} \label{eqdefgenhnwiththeta}
g_\Theta(t)
:=
\sum_{k=1}^\infty\frac{\theta_k}{k} t^k
\quad\mbox{and}\quad
G_\Theta(t)
:=
\exp\Biggl( \sum_{k=1}^\infty\frac{\theta_k}{k} t^k \Biggr).
\end{equation}
At the moment, $g_\Theta(t)$ and $G_\Theta(t)$ are just formal
power series, however, we will see in Section~\ref{secgenfkt}
that
%
\begin{equation}
G_\Theta(t)
=
\sum_{N=0}^\infty h_N t^N,
\label{eqdefgenGtwithhn}
\end{equation}
where the $h_N$ are given in
Definition~\ref{defweightedprobabililtymeasure}.

\subsection{Generating functions and singularity analysis}
\label{secgenfkt}
The idea of generating functions is to encode information of a
sequence into a formal power series.
%
\begin{definition}
\label{defgnerantingfunction} Let $(g_N )_{N\in\N}$
be a sequence of complex numbers and define the (ordinary)
generating function of the sequence as the formal power series
%
\begin{equation}
G(t) = \sum_{N=0}^\infty g_N t^N.
\end{equation}
We define $[t^N][ G ]$ to be the coefficient of $t^N$ in $G(t)$,
that is,~$[t^N][ G ] := g_N$.
\end{definition}

The reason why generating functions are useful is that it is
often possible to compute the generating function without
knowing $g_N$ explicitly.

The main tool in this paper to calculate generating functions
is the following lemma.
%
\begin{lemma}\label{lemcycleindextheorem}
Let $(a_m)_{m\in\N}$ be a sequence of complex numbers. Then
%
\begin{equation}
\label{eqsymmfkt}\qquad
\sum_{\la}\frac{1}{z_\la}\Biggl(\prod_{m=1}^{\ell(\la)} a_{\la_m}\Biggr)
t^{|\la|}
=
\sum_{\la}\frac{1}{z_\la}\Biggl(\prod_{k=1}^{\infty} (a_k t^k)^{C_k} \Biggr)
=
\exp\Biggl(\sum_{k=1}^\infty\frac{1}{k} a_k t^k\Biggr)
\end{equation}
with the same $z_\la$ as in
Lemma~$\ref{lemsizeofconjclasses}$.

If any one of the sums in \eqref{eqsymmfkt} is absolutely
convergent, then so are the others.
\end{lemma}

\begin{pf}
The first equality follows immediately from the definition of
$C_k$. The proof of the second equality in \eqref{eqsymmfkt}
can be found in~\cite{macdonald} or can be directly verified
using the definitions of $z_\la$ and the exponential function.
The last statement follows with dominated convergence.
\end{pf}

We now use this lemma to prove the identity given in
\eqref{eqdefgenGtwithhn}. The constant $h_N$ in~\eqref{eqPThetawithpartition} is chosen so that
$\PT[\sigma]$\vadjust{\goodbreak} is a probability measure on $\SN$. It thus
follows that
%
\begin{equation}
h_N
=
\frac{1}{N!} \sum_{\sigma\in\SN} \prod_{k=1}^{N} \theta_{k}^{C_k}
=
\sum_{\la\vdash N} \frac{1}{z_\la} \prod_{m=1}^{\ell(\la)}
\theta_{\la_m}.
\end{equation}
It now follows, with Lemma~\ref{lemcycleindextheorem}, that
%
\begin{equation} \label{eqexpgenhN}
\sum_{N=0}^\infty h_N t^N
=
\sum_{\la} \frac{1}{z_\la} t^{|\la|} \prod_{m=1}^{\ell(\la)}
\theta_{\la_m}
=
\exp\Biggl( \sum_{k=1}^\infty\frac{\theta_k}{k} t^k \Biggr)
=
G_\Theta(t),
\end{equation}
which proves \eqref{eqdefgenGtwithhn}.

\begin{corollary}\label{corGThetaEwens}
In the special case of the Ewens measure, when $\theta_k$ is a
constant $\theta$, say, we have
%
\begin{equation}
G_\Theta(t) = \sum_{N=0}^\infty h_N t^N = (1-t)^{-\theta} .
\end{equation}
From this it immediately follows that $h_N = (-1)^N
{-\theta\choose N} = {N+\theta-1\choose N}$.
\end{corollary}

Given a generating function $G(t)$, a natural question is: what
is the coefficient of $t^N$ and what is the asymptotic
behavior of this coefficient as $N \to\infty$? If $G(t)$ is
holomorphic near $0$, then one can use Cauchy's integral formula
to do this. But it turns out that it is often difficult to
compute the integral exactly, but we will now see that
one can nevertheless extract the asymptotic behavior of the coefficient
when $G(t)$ has a special form.

\begin{definition} \label{defdelta0}
Given $R>r$ and $0 < \phi<\frac{\pi}{2}$, let
%
\begin{equation}
\Delta_0 = \Delta_0(r,R,\phi) = \{z\in\C\dvtx|z|<R, z \neq r
,|\arg(z-r)|>\phi\}.
\label{eqdefdelta0}
\end{equation}
\end{definition}

\begin{figure}

\includegraphics{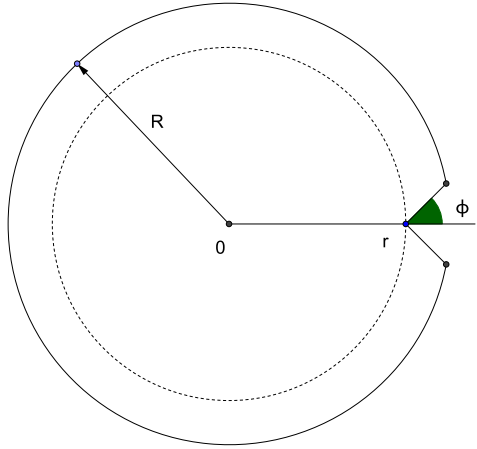}

\caption{Illustration of $\Delta_0(r,R,\phi)$.}\label{figdelta0}
\end{figure}

The domain $\Delta_0$ is
illustrated in Figure~\ref{figdelta0}.

\begin{definition} \label{deffunctionclassFalphar}
Let $r>0$, $\vartheta\geq0$ and a complex constant $K$ be
given. We say that a function $g(t)$ is in $\mathcal{F}(r,\vartheta,
K)$ if there exists $R>r$ and $0 < \phi<\frac{\pi}{2}$ such
that $g(t)$ is holomorphic in $\Delta_0(r,R,\phi)$, and
%
\begin{equation}\label{eqclassFralphanearr}
g(t) = \vartheta\log\biggl( \frac{1}{1-t/r} \biggr) + K + O ( t-r )
\end{equation}
as $t\to r$ with $t\in\Delta_0(r,R,\phi)$.
\end{definition}

%
%

The following theorem, proven by Hwang in
\cite{hwangthesis}, gives the asymptotic behavior of the
coefficient of $t^N$ for certain special generating functions.
%
\begin{theorem}[(Hwang~\cite{hwangthesis})]\label{thmhwang}
Let $g(t) \in\mathcal{F}(r,\vartheta,K)$, and let $S(t)$ be
holomorphic in $|t| \leq r$. Set $G(t,w) = e^{w g(t)} S(t)$,
then
%
\begin{equation} \label{eqthmhwang}
[t^N][ G(t,w) ]
=
\frac{e^{Kw} N^{w\vartheta-1} }{r^N}
\biggl(\frac{S(r)}{\Gamma(\vartheta w)} + O\biggl( \frac{1}{N}\biggr)\biggr)
\end{equation}
uniformly for bounded complex $w$.
\end{theorem}

\begin{remark*}
The idea of the proof is to take a suitable Hankel contour and
to estimate the integral over each piece. The details can be
found in~\cite{hwangthesis}, Chapter 5.
\end{remark*}

\begin{remark*}
One can compute lower order error terms if one has more terms
in the expansion of $g(t)$ near $r$.
\end{remark*}

As a first simple application of this result, we compute the
asymptotic behavior of $h_N$ for the generalized Ewens measure
if $g_\Theta(t)$, as defined in
\eqref{eqdefgenhnwiththeta}, is in
$\mathcal{F}(r,\vartheta,K)$.
%
\begin{lemma}\label{lemasymptotichn}
Let $g_\Theta(t) \in\mathcal{F}(r,\vartheta,K)$. We then have
%
\begin{equation} \label{eqasymptotichn}
h_N = \frac{e^K N^{\vartheta-1}}{r^N }\biggl(\frac{1}{\Gamma(\vartheta)}
+ O\biggl(\frac{1}{N}\biggr)\biggr).
\end{equation}
\end{lemma}

\begin{pf}
We have proven in \eqref{eqexpgenhN} that
$\sum_{N=0}^\infty h_N t^N = \exp( g_\Theta(t))$. We thus can
apply Theorem~\ref{thmhwang} with $g(t) = g_\Theta(t)$, $w=1$
and $S(t) \equiv1$.
\end{pf}

\begin{remark*}
For the Ewens measure, when $\Theta$ is the constant sequence
$(\theta)_{k=1}^\infty$, we have $g_\Theta(t) \in
\mathcal{F}(1,\theta,0)$ and thus $h_N =
\frac{N^{\theta-1}}{\Gamma(\theta)}
(1+O(\frac{1}{N}))$.\vadjust{\goodbreak} However, in this
special case, one can do much more since $h_N$ is known to
equal ${N+\theta-1\choose N}$.
\end{remark*}

%
\begin{figure}

\includegraphics{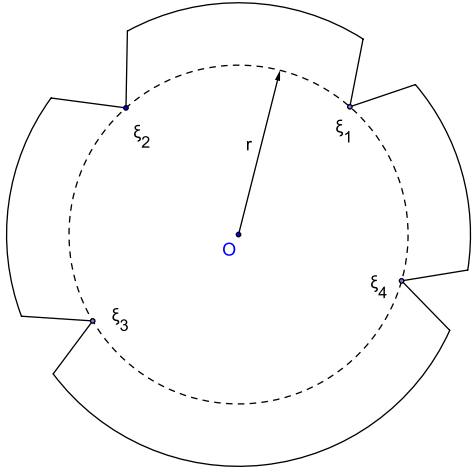}

\caption{Illustration of the domain $\Delta_d(r,R,\phi,\pmb{\xi})$.}\label{figdeltad}
\end{figure}

Essentially, one can think of Hwang's result as concerning
functions with a solitary singularity at $t=r$. In
Section~\ref{secasyptoticmomentsWnf} we will need a
version of this theorem with multiple singularities that we now state.

\begin{definition}\label{defnDeltad}
Let $\pmb{\xi}=(\xi_i)_{i=1}^d$ with $\xi_i \neq\xi_j$ for
$i\neq j$, and with $|\xi_i|=r$ (so the $\xi_i$ are distinct
points lying on the circle of radius $r$). Let $R>r$ and let $0
< \phi<\frac{\pi}{2}$, then set
%
\begin{eqnarray}
&&\Delta_d(r,R,\phi,\pmb{\xi})
\nonumber
\\[-10pt]
\\[-10pt]
\nonumber
&&\qquad:= \bigcap_{i=1}^d \{z\in\C\dvtx
|z|<R, z \neq\xi_i ,|\arg(z-\xi_i)-\arg(\xi_i)|>\phi\}.
\end{eqnarray}
\end{definition}

An example of a
$\Delta_d(r,R,\phi,\pmb{\xi})$ domain is given in Figure~\ref{figdeltad}.

\begin{definition} \label{deffunctionclassFthetaxibr}
Let $\pmb{\vartheta} = (\vartheta_i)_{i=1}^d$ and
$\pmbb{K} = (K_i)_{i=1}^d$ be two sequences of
complex numbers, and let $r>0$. We say a function $g(t)$ is in
$\mathcal{F}(r,\pmb{\vartheta}, \pmbb{K})$ if there exists $R>r$
and $0 < \phi<\frac{\pi}{2}$ such that $g(t)$ is holomorphic
in $\Delta_d(r,R,\phi,\pmb{\xi})$, and for each $i=1,\ldots,d$,\vspace*{-1pt}
%
\begin{equation}\label{eqclassFralphabetanearr}
g(t) = \vartheta_i \log\biggl( \frac{1}{1-t/\xi_i} \biggr) + K_i + O ( t-\xi_i )\vspace*{-1pt}
\end{equation}
as $t\to\xi_i$ with $t\in\Delta_d(r,R,\phi,\pmb{\xi})$.
\end{definition}

Theorem~\ref{thmhwang} generalizes to the next theorem.\vadjust{\goodbreak}
%
\begin{theorem}
\label{thmhwang4} Let\vspace*{1pt} $g\in\mathcal{F}(r,\pmb{\vartheta},
\pmbb{K})$, and let $S(t)$ be holomorphic in $t$ for $|t| \leq
r$. Set $G(t,w) = e^{w g(t)} S(t)$. We have\vspace*{-1pt}
%
\begin{equation} \label{eqthmhwang4}
[t^N][ G(t,w) ]
=
\sum_{i=1}^d \frac{e^{K_i w} N^{w \vartheta_i-1} } {\xi_i^N}
\biggl(\frac{S(\xi_i)}{\Gamma(\vartheta_i w)} + O\biggl( \frac{1 }{N}\biggr)\biggr)
\end{equation}
uniformly for bounded $w$.
\end{theorem}

\begin{pf*}{Sketch of the proof}
The proof is a combination of the proof of a multiple
singularities theorem in~\cite{FlSe09}, Section VI.5, and the
proof of Theorem~\ref{thmhwang}. More precisely, we apply
Cauchy's integral formula with the curve $C$ illustrated in
Figure~\ref{figdeltad+curve}, where the radius $R$ of the
great circle is chosen fix with $R>r$, while the radii of the
small circles are $1/n$.
%
\begin{figure}

\includegraphics{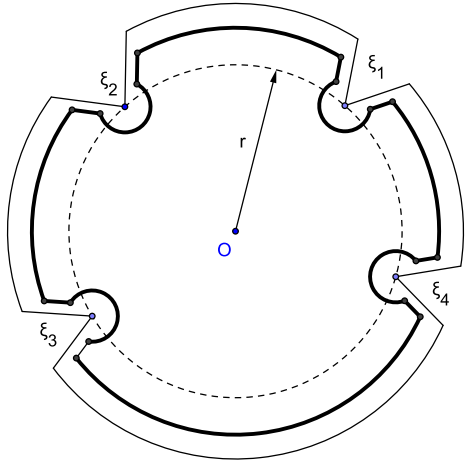}

\caption{The curve $C$.}\label{figdeltad+curve}\vspace*{-3pt}
\end{figure}

A straightforward computation then shows that the integral over
this curve gives~\eqref{eqthmhwang4} and that the error
terms are uniform for bounded $w$.
\end{pf*}

In practice, the computation of the asymptotic behavior near
the singularity is often very difficult, and it is not easy to
prove whether a function $g(t)$ is in
$\mathcal{F}(r,\vartheta,K)$ or not. An alternate approach is
to combine singularity analysis with more elementary methods.
The idea is to write $G = G_1 G_2$ in a way that we can apply
singularity analysis on $G_1$ and can estimate the growth rate
of $[t^N][ G_2 ]$. One then can compute the coefficient $[t^N][ G ]$
directly and apply elementary analysis on it. This method is
called the convolution method.\looseness=-1

\begin{definition}
\label{deffunctionclasseFalphargamma}
Let $\vartheta\geq0, r>0, 0<\gamma\leq1$ be given. We say
$g(t)$ is in $e\mathcal{F}(r,\vartheta,\gamma)$ if $g(t)$ is
holomorphic in $|t|<r$ with
%
\begin{equation}\label{eqclasseFalphargammanearr}
g(t) = \vartheta\log\biggl( \frac{1}{1-t/r} \biggr) + g_0(t)\vadjust{\goodbreak}
\end{equation}
and
%
\begin{equation}
[t^N][ g_0 ] = O(r^{-N} N^{-1-\gamma})
\end{equation}
as $N\to\infty$.
\end{definition}
%
\begin{theorem}[(Hwang~\cite{MR1724562})]\label{thmhwang2}
Let $g(t) \in e\mathcal{F}(r,\vartheta,\gamma)$, and let $S(t)$
be holomorphic in $|t| \leq r$. Set $G(t,w) = e^{w g(t)} S(t)$,
then
%
\begin{equation}
[t^N][ G(t,w) ]
=
\frac{e^{w g_0(r)} N^{w\vartheta-1} }{r^N} \frac{S(r)}{\Gamma
(\vartheta w)} + R_N(w)
\end{equation}
with
%
\begin{equation}
R_N(w)
=
\cases{
\displaystyle O\biggl(\frac{N^{\vartheta\Re(w) - 1- \gamma} \log(N)}{r^N} \biggr), &\quad
$\mbox{if } \Re(w) \geq0,$
\vspace*{2pt}\cr
\displaystyle O\biggl(\frac{ N^{-1 - \gamma} }{r^N} \biggr), &\quad$\mbox{if } \Re(w) < 0,$
}
\end{equation}
uniformly for bounded $w$.
\end{theorem}

This theorem is more general than Theorem~\ref{thmhwang}, but
the error terms are worse. As in Lemma~\ref{lemasymptotichn},
we can compute the asymptotic behavior of $h_N$ in the case of the
generalized Ewens measures when $g_\Theta(t) \in
e\mathcal{F}(r,\vartheta,\gamma)$.
%
\begin{lemma} \label{lemasymptotichn2}
Assume that $g_\Theta(t) \in e\mathcal{F}(r,\vartheta,\gamma)$. We
then have
%
\begin{equation}
h_N = \frac{e^{g_0(r)} N^{\vartheta-1}}{r^N \Gamma(\vartheta)} +
O\biggl(\frac{N^{\vartheta-1-\gamma} \log N}{r^N}\biggr).
\end{equation}
\end{lemma}

\section{Moments of multiplicative class functions}
\label{sectACF}
We extend in this section the results of~\cite{HKOS,Zei10} and
\cite{Zei10a} to the generalized Ewens measure $\PT$. More
precisely, we compute the asymptotic behavior of the moments
of the characteristic polynomial $Z_N(x)$ and of multiplicative
class functions $W^N(P)$ with respect to $\PT$ using the methods
of generating functions and singularity analysis introduced in
the previous section.

\subsection{Multiplicative class functions}
It is well known that $\SN$ can be identified with the group of
permutation matrices via
%
\begin{equation}\label{eqembddingSn}
\sigma\mapsto\bigl(\delta_{i,\sigma(j)}\bigr)_{1\leq i,j\leq N}.
\end{equation}
It is easy to see that this map is an injective group
homomorphism. We thus do not distinguish between $\SN$ and the
group of permutation matrices and use for both the notation
$\SN$. It will always be clear from the context if it is
necessary to consider $\sigma\in\SN$ as a matrix.

\begin{definition}
Let $x\in\C$ and $\sigma\in\SN$. The characteristic polynomial
of $\sigma$ is
%
\begin{equation}\label{eqdefzndet}
Z_N(x)=Z_N(x)(\sigma):=\det(I_N-x\sigma).
\end{equation}
\end{definition}

It is a standard fact that the characteristic polynomial can be written
in terms of the
cycle type of $\sigma$.
%
\begin{lemma}
Let $\sigma\in\SN$ be given with cycle type $\la$; then
%
\begin{equation}
\label{eqdefznla}
Z_N(x) = \prod_{m=1}^{\ell(\la)} (1-x^{\la_m}),
\end{equation}
with $\ell(\la)$ the length of the partition $\lambda$, which
is the same as the number of cycles~$T(\sigma)$.
\end{lemma}
\begin{pf}
Since any permutation matrix is conjugate to a block matrix with
each block corresponding to one of the cycles, and the
characteristic polynomial factors over the blocks, it is
sufficient to prove this result in the simple case of a
one-cycle permutation, where it follows from a simple
calculation. More explicit details can be found, for instance,
in~\cite{Zei10a}, Chapter~2.2.
\end{pf}

Equation~\eqref{eqdefznla} shows that the spectrum of
permutation matrix is uniquely determined by the cycle type. We
use this as motivation to define multiplicative class functions
on $\SN$.
%
\begin{definition}
\label{defWpolyinone}
Let $P(x)$ be a polynomial in $x$. We then define the
multiplicative class function associated to the polynomial $P$
as
%
%
\begin{equation}
\label{eqdefWpolyinone}
W^N(P)(x)
=
W^N(P)(x)(\sigma)
:=
\prod_{m=1}^{\ell(\la)}P(x^{\la_m}).
\end{equation}
For brevity, we simply call this a multiplicative class
function.
\end{definition}

It follows immediately that the characteristic polynomial is
the multiplicative class function associated to the polynomial
$P(x) = 1-x$. The main difference between $Z_N(x)$ and $W^N(P)$
is that $W^N(P)$ is independent of the interpretation of $\SN$
as matrices.

We now wish to obtain the asymptotic behavior of
the moments
\[
\mathbb{E}_{\Theta}[( W^N(P_1)(x_1) )^{k_1} ( W^N(P_2)(x_2) )^{k_2}]
\]
for $x_1 \neq x_2$. The easiest
way to achieve this is to extend the definition of $W^{N}(P)$.
%
\begin{definition}
\label{defW1poly}
Let $P(x_1,x_2)$ be a polynomial in the two variables
$x_1,x_2$. For $\sigma\in\SN$ with cycle type $\la$, we set
%
\begin{equation}\label{eqdefWpolyintwo}
W^N(P)(x_1,x_2)
=
W^N(P)(x_1,x_2)(\sigma)
:=
\prod_{m=1}^{\ell(\la)}P(x_1^{\la_m}, x_2^{\la_m}).
\end{equation}
\end{definition}

A simple computation using the definitions above shows
\begin{eqnarray*}
( W^N(P_1)(x_1) )^{k_1} ( W^N(P_2)(x_2) )^{k_2}
&=&
W^N(P^{k_1}_1)(x_1) W^N(P^{k_2}_2)(x_2)\\
&=&
W^N(P)(x_1,x_2)
\end{eqnarray*}
with $P(x_1,x_2) = P^{k_1}_1(x_1) P^{k_2}_2(x_2)$.

This shows that it is enough to consider
$\mathbb{E}_{\Theta}[W^N(P)(x_1,x_2)]$.

\begin{remark*}
There is no restriction to the number of variables. All
arguments used here work also for more than two variables. We
restrict ourselves to two variables since this is enough to
illustrate the general case.
\end{remark*}

\subsection{Generating functions for $W^N$}
\label{secgeneratingforasso}
In this section we compute the generating functions of the
moments of multiplicative class functions.

%
\begin{theorem}
\label{thmgeneratingassociated}
Let $P$ be a complex polynomial with
%
\begin{equation}
P(x_1,x_2)
=
\sum_{k_1,k_2=0}^{\infty}b_{k_1,k_2} x_1^{k_1} x_2^{k_2}.
\end{equation}
We have as formal power series
%
\begin{equation}
\sum_{N=0}^\infty t^N h_N \mathbb{E}_{\Theta}[W^{N}(P)(x_1,x_2)]
=
\prod_{k_1,k_2=0}^{\infty} ( G_\Theta(x_1^{k_1} x_2^{k_2} t)
)^{b_{k_1,k_2} },
\label{eqgeneratingW1poly}
\end{equation}
where $G_\Theta(t)^b = \exp(b \cdot g_\Theta(t))$, and
$G_\Theta$ and $g_\Theta$ are defined in
\eqref{eqdefgenhnwiththeta}.
\end{theorem}

\begin{pf}
From \eqref{eqdefWpolyintwo} and
\eqref{eqPThetawithpartition} we have
%
\begin{equation}
\mathbb{E}_{\Theta}[W^{N}(P)(x_1,x_2)]
= \frac{1}{h_N N!} \sum_{\sigma\in\SN}\prod_{m=1}^{\ell(\la)}
\theta_{\la_m} P(x_1^{\la_m},x_2^{\la_m})
\end{equation}
and since $W^N(P)$ is a class function, we may use
\eqref{eqsumclassfuncs} to obtain
%
\begin{equation}
\mathbb{E}_{\Theta}[W^{N}(P)(x_1,x_2)] = \frac{1}{h_N} \sum_{\la
\vdash N}\frac{1}{z_\la}\prod_{m=1}^{\ell(\la)} \theta_{\la_m}
P(x_1^{\la_m},x_2^{\la_m}).
\end{equation}
We now compute the generating function of $h_N \mathbb{E}_{\Theta
}[W^{N}(P)(x_1,x_2)]$ with the help of
Lemma~\ref{lemcycleindextheorem}:
\begin{eqnarray*}
\sum_{N=0}^\infty t^{N} h_N\mathbb{E}[W^{N}(P)(x_1,x_2)] &=& \sum
_{N=0}^\infty t^N \sum_{\la\vdash N}\frac{1}{z_\la}\prod
_{m=1}^{\ell(\la)} \theta_{\la_m} P(x_1^{\la_m},x_2^{\la_m}) \\
&= &\sum_{\la}\frac{1}{z_\la}t^{|\la|} \prod_{m=1}^{\ell(\la)}
(\theta_{\la_m} P(x_1^{\la_m},x_2^{\la_m}) )\\
&=&
\exp\Biggl( \sum_{m=1}^\infty\frac{\theta_m}{m} t^m P(x_1^{m},x_2^{m})
\Biggr)\\
&=&
\exp\Biggl( \sum_{k_1,k_2=0}^{\infty} b_{k_1,k_2} \sum_{m=1}^\infty\frac
{\theta_m}{m} (x_1^{k_1} x_2^{k_2} t)^m \Biggr).
\end{eqnarray*}
Note that
%
\begin{equation}
\sum_{m=1}^\infty\frac{\theta_m}{m} (x_1^{k_1} x_2^{k_2} t)^m =
g_\Theta( x_1^{k_1} x_2^{k_2} t ),
\end{equation}
where $g_\Theta$ is defined in
\eqref{eqdefgenhnwiththeta}, and hence,
%
\begin{equation}
\sum_{N=0}^\infty t^{N} h_N\mathbb{E}[W^{N}(P)(x_1,x_2)]
= \prod_{k_1,k_2=0}^{\infty} (
G_\Theta(x_1^{k_1} x_2^{k_2} t) )^{b_{k_1,k_2} }
\end{equation}
and this proves \eqref{eqgeneratingW1poly}, as required.
\end{pf}

\begin{remark*}
The requirement for $P$ to be a polynomial is there to ensure
absolute convergence, and clearly this condition can be
considerably weakened (see~\cite{Zei10}, Sections~5, 6).
\end{remark*}

As an immediate consequence we have the following corollary.
%
\begin{corollary}
\label{corgeneratingZn} Let $s_1,s_2\in\N$ be given. Then
%
\begin{eqnarray}\label{eqgeneratingZn}
&&\sum_{N=0}^\infty t^N h_N \mathbb{E}_{\Theta}[Z_N(x_1)^{s_1} Z_N(x_2)^{s_2}]
\nonumber
\\[-8pt]
\\[-8pt]
\nonumber
&&\qquad
=
\prod_{k_1=0}^{s_1} \prod_{k_2=0}^{s_2} ( G_\Theta(x_1^{k_1}
x_2^{k_2} t) )^{(-1)^{k_1+k_2} {s_1 \choose
k_1}{s_2\choose
k_2}}.
\end{eqnarray}

\end{corollary}
\begin{pf}
We have
\begin{eqnarray*}
Z_N(x_1)^{s_1} Z_N(x_2)^{s_2}
&=&
\bigl( W^{N}(1-x_1) \bigr)^{s_1}\bigl( W^{N}(1-x_2) \bigr)^{s_2}\\
&=&
W^{N} \bigl( (1-x_1)^{s_1} (1-x_2)^{s_2} \bigr).
\end{eqnarray*}
The corollary now follows immediately by calculating the Taylor
expansion of $(1-x_1)^{s_1} (1-x_2)^{s_2}$ near $0$.
\end{pf}
%

\subsection{Asymptotic behavior of the moments}
\label{secasyptoticmomentsWnf}
Combining the generating functions in
Theorem~\ref{thmgeneratingassociated} with the singularity
analysis developed in Section~\ref{secgenfkt}, we compute the
asymptotic behavior of $\mathbb{E}_{\Theta}[W^{N}(P)]$ as $N\to
\infty$.

We have to distinguish between the cases $|x_i|<1$ and $|x_i|=1$. We
consider here only $g_\Theta(t)\in\mathcal{F}(r, \vartheta,K)$.
The results and computations for $g_\Theta(t)\in
e\mathcal{F}(r, \vartheta,\gamma)$ are similar, with only minor
differences in the error terms.

We first look at the asymptotic behavior inside the unit disc.
We have the following theorem.
%
\begin{theorem}
\label{thmasymptoticmomentsassociated}
Let $P$ be as in Theorem~$\ref{thmgeneratingassociated}$, and
let $x_1,x_2\in\C$ be given with $\max\{|x_1|,|x_2|\} <
1$.
Assume that $g_\Theta(t) \in\mathcal{F}(r,\vartheta,K)$, then
%
\begin{equation}\label{eqasymptoticinsideunitW1}
\mathbb{E}[W^{N}(P)]
=
N^{\vartheta(b_{0,0}-1)} e^{K(b_{0,0}-1)}
\biggl(
E_1 + O\biggl(\frac{1}{N}\biggr) \biggr),
\end{equation}
with
%
\begin{equation}
E_1
=
E_1(x_1,x_2)
=
\frac{\Gamma(\vartheta)}{\Gamma(\vartheta b_{0,0})}
\prod_{(k_1,k_2) \neq(0,0)} ( G_\Theta(r x_1^{k_1} x_2^{k_2})
)^{b_{k_1,k_2}}.
\end{equation}
\end{theorem}

\begin{pf}
Set
%
\begin{equation}
S(t):= \prod_{(k_1,k_2) \neq(0,0)} ( G_\Theta(x_1^{k_1} x_2^{k_2} t)
)^{b_{k_1,k_2} }.
\label{eqStinsidethecircle}
\end{equation}
Since $P$ is polynomial, the product is finite and there is no
problem with convergence. The domain of holomorphicity of $S$
is thus the intersection of the domains of holomorphicity of
each factor. This shows that the function $S(t)$ is holomorphic
for $|t|< r +\varepsilon$ for an $\varepsilon>0$ since
$\max\{|x_1|,|x_2|\}<1$ and $G_\Theta(t)$ is holomorphic for
$|t|<r$.

Separating the $k_1=k_2=0$ term in
\eqref{eqgeneratingW1poly} from the rest, we can write the
generating function as
%
\begin{equation}
\sum_{N=0}^\infty t^N h_N \mathbb{E}_{\Theta}[W^{N}(P)(x_1,x_2)]
= \exp\bigl( b_{0,0} \cdot g_\Theta(t) \bigr) S(t).
\end{equation}

Applying Theorem~\ref{thmhwang}, we get
\[
h_N \mathbb{E}[W^{N}(P)]
=
N^{\vartheta b_{0,0}-1} e^{K b_{0,0}} \frac{1}{r^N} \biggl( \frac
{S(r)}{\Gamma(\vartheta b_{0,0})} + O\biggl(\frac{1}{N}\biggr) \biggr).
\]
Comparing $S(r)$ with $E_1$, and using
Lemma~\ref{lemasymptotichn} to find the asymptotic behavior
of $h_N$, proves the theorem.
\end{pf}

As a special case, we get the asymptotic behavior of
$\mathbb{E}_{\Theta}[Z_N^{s_1}(x_1) Z_N^{s_2}(x_2) ]$ with respect to
$\PT$
inside the unit disc.

\begin{corollary}
\label{corZNinside}
Let $x_1,x_2\in\C$ be given with $\max\{|x_1|,|x_2|\} < 1$ and let
$s_1,s_2\in\N$. We then have
%
\begin{eqnarray}
&&\mathbb{E}_{\Theta}[Z_N^{s_1}(x_1) Z_N^{s_2}(x_2) ]
\nonumber
\\[-8pt]
\\[-8pt]
\nonumber
&&\qquad=
\prod_{(k_1,k_2) \neq(0,0)} ( G_\Theta(r x_1^{k_1} x_2^{k_2})
)^{(-1)^{k_1+k_2} {s_1\choose k_1} {s_2\choose k_2}}
+ O\biggl(\frac{1}{N}\biggr).
\end{eqnarray}
\end{corollary}

\begin{pf}
This follows immediately from the fact that $Z_N(x) =W^N(1-x)(x)$ and
that $(1-x_1)^{s_1} (1-x_2)^{s_2}$ evaluated at $x_1 =x_2=0$ is $1$.
\end{pf}

In particular, for the uniform measure ($\theta_k \equiv1$ for
all $k$) Corollary~\ref{corGThetaEwens} gives $G_\Theta(t) =
(1-t)^{-1}$ in which case we have
%
\begin{eqnarray}
&&\mathbb{E}[Z_N^{s_1}(x_1) Z_N^{s_2}(x_2) ]
\nonumber
\\[-8pt]
\\[-8pt]
\nonumber
&&\qquad=
\prod_{(k_1,k_2) \neq(0,0)} (1-x_1^{k_1} x_2^{k_2} )^{-
{s_1\choose k_1} {{s_2\choose k_2}}(-1)^{k_1+k_2}} + O\biggl(\frac{1}{N}\biggr).
\end{eqnarray}
This shows that Corollary~\ref{corZNinside} agrees with
\cite{Zei10a}, Theorem 2.13, in the unform case.

The behavior on the unit disc is more complicated. The reason
is that the generating function can have (for fixed $x_1,x_2$)
more than one singularity on the circle of radius $r$. Another
point that makes this case more laborious is the requirement to
check whether some of the singularities of the factors on the
right-hand side of~\eqref{eqgeneratingW1poly} are equal. For simplicity,
we assume that all singularities are distinct.
%
\begin{theorem}
\label{thmasymptoticonunitdiscW1} Let $P$ be as in
Theorem~$\ref{thmgeneratingassociated}$, and let $x_1,x_2\in\C
$ be given with $|x_1|= |x_2|=1$ and $x_1^{k_1} x_2^{k_2} \neq
1$ for all $(k_1,k_2)\in\Z^2\setminus\{(0,0)\}$.
Assume that $g_\Theta(t) \in\mathcal{F}(r,\vartheta,K)$, then
\begin{eqnarray*}
&&\mathbb{E}_{\Theta}[W^{N}(P)] \\
&&\qquad=
\mathop{\sum_{k_1,k_2}}_{b_{k_1,k_2} \neq0}
E_2(k_1,k_2)
N^{\vartheta(b_{k_1,k_2} -1)} x_1^{N k_1} x_2^{N k_2}
\biggl(\frac{\Gamma(\vartheta)}{\Gamma(\vartheta b_{k_1,k_2} )}
+ O\biggl(\frac{1}{N}\biggr) \biggr)
\end{eqnarray*}
with
%
\begin{equation}
 E_2(k_1,k_2)
=
e^{K (b_{k_1,k_2} -1)}
\prod_{(m_1,m_2) \neq(k_1,k_2)} ( G_\Theta(rx_1^{m_1-k_1}
x_2^{m_2-k_2}) )^{b_{m_1,m_2}}.\hspace*{-35pt}
\end{equation}
\end{theorem}

\begin{pf}
We define
%
\begin{equation}
F(t)
=
\sum_{k_1,k_2} b_{k_1,k_2} g_\Theta(x_1^{k_1} x_2^{k_2} t) .
\end{equation}
By \eqref{eqgeneratingW1poly} we see that $\exp(F(t))$ is
the generating function of $h_N \mathbb{E}_{\Theta}[W^{N}(P)]$. We
first take
a look at the domain of holomorphicity of $F(t)$. We have by
assumption that $g_\Theta(t)$ is holomorphic in
$\Delta_0(r,R,\phi)$ for an $R>r$ and $0 < \phi<\frac{\pi}{2}$.
This shows that $g_\Theta(x_1^{k_1} x_2^{k_2} t)$ is
holomorphic for $t$ in the domain $\Delta_1(r,R,\phi,r
x_1^{-k_1} x_2^{-k_2})$ with $\Delta_1$ as in
Definition~\ref{defnDeltad}, and that $F$ is holomorphic in
%
\begin{equation}
D:=
\mathop{\bigcap_{k_1,k_2}}_{ b_{k_1,k_2} \neq0} \Delta_1(r,R,\phi,
r x_1^{-k_1} x_2^{-k_2})
=
\Delta_d(r,R,\phi,\pmb{\xi}),
\end{equation}
where $\pmb{\xi}$ is the finite sequence of all $rx_1^{-k_1}
x_2^{-k_2}$ with $b_{k_1,k_2} \neq0$ (in any order).
Notice that this is only a finite intersection since $P$ is a
polynomial. Since $|x_1| = |x_2| =1$, we see that $D$ has a
shape as in Figure~\ref{figdeltad} and that $F$ has
singularities at $t= r x_1^{-k_1} x_2^{-k_2}$. We thus may use
Theorem~\ref{thmhwang4} and therefore need to take a look at
the behavior of $F$ near each singularity. We assumed that
$x_1^{k_1} x_2^{k_2} \neq x_1^{m_1} x_2^{m_2}$ for $(m_1,m_2)
\neq(k_1,k_2)$, which implies that the singularities are
distinct, and thus near the point $r x_1^{-k_1} x_2^{-k_2}$,
$F(t)$ has the expansion
%
\begin{eqnarray}\label{eqexpansionFnearsing}
\qquad F(t)
&=&
b_{k_1,k_2} \vartheta\log\biggl( \frac{1}{1-tx_1^{k_1} x_2^{k_2}/r} \biggr)
+
b_{k_1,k_2}K
\nonumber
\\[-8pt]
\\[-8pt]
\nonumber
&&{}+
\sum_{(m_1,m_2) \neq(k_1,k_2)} b_{m_1,m_2} g_\Theta( r x_1^{m_1-k_1}
x_2^{m_2-k_2} ) + O(t-r x_1^{-k_1} x_2^{-k_2})
\end{eqnarray}
for $t \to r x_1^{-k_1} x_2^{-k_2}$. This shows that we can
apply Theorem~\ref{thmhwang4}. Combining this together with
Lemma~\ref{lemasymptotichn} proves the theorem.
\end{pf}

\begin{remark*}
For simplicity we have assumed that all the singularities are
distinct. The modification required to cope with the case when
$x_1^{k_1} x_2^{k_2} = x_1^{m_1} x_2^{m_2}$ for some $(m_1,m_2)
\neq(k_1,k_2)$ would appear in
\eqref{eqexpansionFnearsing}, but technically there is no
restriction. Such a situation, with all the details written out
explicitly, appears in~\cite{Zei10}.
\end{remark*}

To illustrate this theorem, we will calculate the
autocorrelation of two characteristic polynomials at distinct
points $x_1,x_2$ on the unit circle subject to $x_1^{k_1} \neq
x_2^{k_2}$ for all $\{k_1,k_2\} \neq\{0,0\}$.

The four coefficients of $Z_N(x_1) Z_N(x_2)$ are easy to
calculate, being $b_{0,0} = b_{1,1} = 1$ and $b_{1,0} = b_{0,1}
= -1$ and this enables an immediate simplification to occur by
observing that only the terms with $b_{k_1,k_2}$ maximal
contribute; the others are of lower order, in this case being
$O(N^{-2\vartheta})$.

Substituting these values into the theorem we have
%
\begin{eqnarray}
&&\mathbb{E}_{\Theta}[Z_N(x_1) Z_N(x_2)]
\nonumber
\\[-8pt]
\\[-8pt]
\nonumber
&&\qquad
= E_2(0,0)+
E_2(1,1) x_1^N x_2^N
+O\biggl(\frac{1}{N}\biggr) +O\biggl(\frac{1}{N^{2\vartheta}}\biggr)
\end{eqnarray}
with
%
\begin{equation}
E_2(0,0) = \frac{ G_\Theta(r x_1 x_2 ) }{ G_\Theta(rx_1 ) G_\Theta
(rx_2 ) }
\end{equation}
and
%
\begin{equation}
E_2(1,1) = \frac{ G_\Theta(r x_1^{-1} x_2^{-1} ) }{ G_\Theta
(rx_1^{-1} ) G_\Theta(rx_2^{-1} ) }.
\end{equation}

\section{Traces}\label{sectraces}

\subsection{Traces of permutation matrices}
In this section we consider the asymptotic behavior of traces
of permutation matrices. Powers of traces and traces of powers
have received much attention in the random matrix literature
(see, e.g.,~\cite{DiaEva01,DiaGam04,DiacShah94}). More
specifically, we
first look at $\tr(\sigma^d)$ for fixed $d\in\Z$. Since the
embedding of $\SN$ into the unitary group in
\eqref{eqembddingSn} is a group homomorphism, we can
interpret $\sigma^d$ as $d$-fold matrix multiplication and as
the matrix corresponding to $\underbrace{\sigma\circ\cdots
\circ\sigma}_{d\ \mathrm{times}}$.

We first recall a well-known explicit expression for $\tr(\sigma^d)$
that we shortly prove for completeness.

\begin{lemma}
\label{lemtracespermutationexplicit} We have for $d\in\Z$
%
\begin{equation}\label{eqTracePower}
\tr(\sigma^d) = \sum_{k = 1}^N \one_{k | d} k C_k(\sigma),
\qquad \mbox{with }
\one_{k | d} =
\cases{
1, &\quad$\mbox{if } k \mbox{ divides }d,$\vspace*{2pt}\cr
0, &\quad$\mbox{otherwise}.$}
\end{equation}
\end{lemma}

\begin{pf}
The matrix corresponding to $\sigma^d$ has the form
$(\delta_{i,\sigma^d(j)})$. We thus have
%
\begin{equation}
\tr(\sigma^d) = \sum_{i = 1}^N \delta_{i,\sigma^d(i)} = \#\{i
\dvtx\sigma^d(i) = i\}.
\end{equation}
Therefore, $\tr(\sigma^d)$ is the number of $1$-cycles of
$\sigma^d$. A simple computation now shows that the number of
$1$-cycles of $\sigma^d$ is indeed $\sum_{k = 1}^N \one_{k | d}
k C_k(\sigma)$.
\end{pf}

Using this expression and the method of generating functions
developed in Section~\ref{secgenfkt}, we prove a weak
convergence result for $\tr(\sigma^d)$.
%
\begin{theorem}\label{thmcharfunctracematrix}
Let $d\in\N$ be given. We then have
%
\begin{equation}\label{eqcharfunctracematrix1}
\sum_{N=0}^\infty h_N\mathbb{E}_{\Theta}\bigl[e^{i s\tr(\sigma^d)}\bigr] t^N
=
\exp\biggl(\sum_{k|d} \frac{\theta_k}{k} (e^{i sk} - 1) t^k \biggr) G_\Theta(t).
\end{equation}
If $g_\theta$ is of class $\mathcal{F}(\vartheta,r,K)$, then
%
\begin{equation}\label{eqcharfunctracematrix2}
\mathbb{E}_{\Theta}\bigl[e^{i s\tr(\sigma^d)}\bigr]
=
\exp\biggl( \sum_{k|d} \frac{\theta_k}{k} (e^{i sk}-1) r^k \biggr) + O\biggl( \frac
{1}{N} \biggr).
\end{equation}
If $g_\theta$ is of class $e\mathcal{F}(\vartheta,r,\gamma)$,
then
%
\begin{equation}\label{eqcharfunctracematrix3}
\mathbb{E}_{\Theta}\bigl[e^{i s\tr(\sigma^d)}\bigr]
=
\exp\biggl( \sum_{k|d} \frac{\theta_k}{k} (e^{i sk}-1) r^k \biggr) + O\biggl( \frac
{\log(N)}{N^\gamma} \biggr).
\end{equation}
\end{theorem}

\begin{pf}
Applying Lemma~\ref{lemtracespermutationexplicit}, and
evaluating the expectation explicitly in terms of partitions
using Lemma~\ref{lemsizeofconjclasses}, we have
%
\begin{equation}
\sum_{N=0}^\infty h_N\mathbb{E}_{\Theta}\bigl[e^{i s\tr(\sigma^d)}\bigr] t^N
=
\sum_{\la}\frac{1}{z_\la}\prod_{k=1}^\infty(\theta_k e^{i s k
\one_{k|d} } t^k )^{C_k}.
\end{equation}
The cycle index theorem (Lemma~\ref{lemcycleindextheorem})
yields that this equals
%
\begin{equation}
\exp\Biggl( \sum_{k=1}^\infty\frac{1}{k} \theta_k e^{i s k \one_{k|d}}
t^k \Biggr)
= \exp\biggl( \sum_{k|d} \frac{\theta_k}{k} (e^{i sk} - 1) t^k \biggr) G_\Theta(t),
\end{equation}
where $G_\Theta(t)$ is given in
\eqref{eqdefgenhnwiththeta}. This proves
equation~\eqref{eqcharfunctracematrix1}.

Applying Theorem~\ref{thmhwang} to this yields
equation~\eqref{eqcharfunctracematrix2}, and
Theorem~\ref{thmhwang2} yields
equation~\eqref{eqcharfunctracematrix3}, as required.
\end{pf}

\begin{remark*}
An alternative way to prove
Theorem~$\ref{thmcharfunctracematrix}$ is to use Theorem~3.1
in~\cite{NikZei}, which computes the generating function of\break
$h_N \mathbb{E}_{\Theta}[\exp( \sum_{k=1}^b i s_k C_k )]$ and its
asymptotic behavior for $g_\Theta(t) \in
\mathcal{F}(\vartheta,r,K)$ and $g_\Theta(t) \in
e\mathcal{F}(\vartheta,r,\gamma)$.
\end{remark*}

We obtain the following as an immediate corollary.
%
\begin{corollary}
\label{corasymptotictracesofmatrix} Let $d \in\Z$ be fixed
and assume that $g_\Theta$ is in $\mathcal{F}(\vartheta,r,K)$
or $e\mathcal{F}(\vartheta,r,\gamma)$. Then
%
\begin{equation}
\tr(\sigma^d) \stackrel{d}{\longrightarrow} \sum_{k|d} k P_k \qquad\mbox{as }N\to\infty,
\end{equation}
where $P_k$ are independent Poisson distributed random
variables with\break $\mathbb{E}[P_k]= \frac{\theta_k}{k}r^k$.
\end{corollary}

\subsection{Traces of functions}
Recall from the \hyperref[sec1]{Introduction} that if $M$ is the permutation
matrix representing the permutation $\sigma$, then for a
function $F\dvtx S^1 \to\C$, we defined the trace of $F$ to be
the function $\tr(F)\dvtx\SN\to\C$ with
%
\begin{equation}
\tr(F)(\sigma):= \sum_{k=1}^N F(\omega_k),
\end{equation}
where $(\omega_k)_{k=1}^N$ are the eigenvalues of $M$ or $\sigma$ with
multiplicity. Lemma~\ref{lemtracefunctionexplicit} showed that $\tr(F)$
could be expressed in terms of the cycle structure of $\sigma$ as
%
\begin{equation}\label{eqtrFexplicit}
\tr(F)(\sigma)
=
\sum_{k=1}^{N} kC_k \Delta_{k}(F)
\end{equation}
with
%
\begin{equation}\label{eqDeltak}
\Delta_{k}(F):= \frac{1}{k} \sum_{m=1}^k F( e^{2\pi i m/k}
).
\end{equation}

The asymptotic behavior of $\tr(F)$ is not so easy to compute
for an arbitrary function defined on the unit circle. This
problem will be dealt with more carefully in
Sections~\ref{secgeneralcase} and~\ref{sectEwensDivVar}.
However, if $F$ is a Laurent polynomial, we can use the same
method as for $\tr(\sigma^d)$.

\begin{theorem}
Let
%
\begin{equation}
F(x) = \sum_{d} b_d x^d
\end{equation}
be a Laurent polynomial. If $g_\Theta\in
\mathcal{F}(\vartheta,r,K)$ or $g_\Theta\in
e\mathcal{F}(\vartheta,r,\gamma)$, then
%
\begin{equation} \label{eqcharLaurent}
\tr(F)(\sigma) - N b_0 \stackrel{d}{\longrightarrow} \mathop{\sum
_{d=-\infty}}_{ d\neq0}^\infty b_d
\mathop{\sum_{k \geq1}}_{ k|d} k P_k\qquad \mbox{as }N\to\infty,
\end{equation}
where $P_k$ are independent Poisson distributed random
variables with\break $\mathbb{E}[P_k]= \frac{\theta_k}{k}r^k$.
\end{theorem}

\begin{pf}
Due to the linearity of $\tr(F)$, we may assume the constant
term, $b_0$, is zero. As in the previous computations, we apply the cycle index
theorem to obtain
%
\begin{eqnarray}
\sum_{N=0}^\infty h_N\mathbb{E}_{\Theta}\bigl[e^{i s\tr(F)(\sigma)}\bigr] t^N
&=&
\sum_{\la}\frac{1}{z_\la}\prod_{k=1}^\infty\bigl(\theta_k e^{i s k
\Delta_k(F) } t^k \bigr)^{C_k} \\
&=& \exp\Biggl( \sum_{k=1}^\infty\frac{\theta_k}{k} e^{i s k \Delta
_k(F)} t^k \Biggr)
\end{eqnarray}
and since $\Delta_k(x^d) = \one_{k|d}$ and is linear, this
equals
%
\begin{eqnarray}
&&\exp\Biggl( \sum_{k=1}^\infty\frac{\theta_k}{k} \exp\biggl( i s
\sum_{d\neq0} b_d k \one_{k|d}\biggr) t^k \Biggr)
\nonumber
\\[-8pt]
\\[-8pt]
\nonumber
&&\qquad= \exp\Biggl( \sum_{k=1}^\infty\frac{\theta_k}{k}
\biggl(\mathop{\prod_{d\neq0 }}_{ k|d} e^{i s b_ d k}
- 1 \biggr)t^k \Biggr) G_\Theta(t).
\end{eqnarray}
Note the first factor is entire, so Theorem~\ref{thmhwang}
[for the case of $\mathcal{F}(r,\vartheta,K)$] and
Theorem~\ref{thmhwang2} [for the case of
$e\mathcal{F}(r,\vartheta,\gamma)$] yields
%
\begin{equation}
\mathbb{E}_{\Theta}\bigl[e^{i s\tr(F)(\sigma)}\bigr] \to\exp\Biggl( \sum
_{k=1}^\infty\frac{\theta_k}{k}
\biggl(\mathop{\prod_{d\neq0}}_{ k|d} e^{i s b_ d k}
- 1 \biggr) r^k \Biggr).
\end{equation}
The right-hand side is the characteristic function of the right-hand
side in
\eqref{eqcharLaurent}. The proof is complete.
\end{pf}

%
\begin{theorem}
Let $F\dvtx S^1 \to\C$ be of bounded variation, and $d\in\N$ be
given. Then
%
\begin{equation}
\frac{1}{N^d}\mathbb{E}_{\Theta}[(\tr(F)(\sigma))^d] = \biggl(\int
_{S^1} F(\varphi) \, d\varphi\biggr)^d + O\biggl(\frac{\mathbb{E}_{\Theta
}[T(\sigma)]}{N}\biggr),
\label{eqweakconvSn}
\end{equation}
where $T(\sigma)$ is the total number of cycles of $\sigma$,
and $d\varphi$ the uniform measure on ${S^1}$.

Moreover, if $g_\Theta(t)
\in\mathcal{F}(\vartheta,r,K)$ or $g_\Theta(t) \in
e\mathcal{F}(\vartheta,r,\gamma)$, then $\mathbb{E}_{\Theta
}[T(\sigma)] \sim
\vartheta\log(N)$, and thus we have a quick convergence of
the moments.
\end{theorem}

\begin{pf}
Since $F$ is of bounded variation we can apply Koksma's
inequality
(\cite{kuipers-niederreiter-74}, Theorem 5.1) to see
that
%
%
\begin{equation}
\Biggl| \frac{1}{k} \sum_{m=1}^k F( e^{2\pi i m/k })
-
\int_{{S^1}} F(\varphi) \, d\varphi
\Biggr|
\leq
2D_k V(F)
\end{equation}
with $V(F)$ the variation of $F$ and $D_k $ the discrepancy of
the sequence $( e^{2\pi i m/k} )_{m=1}^k$. But the
discrepancy $D_k$ is dominated by $1/k$. We thus have
%
\begin{equation}
\Delta_k(F) = \frac{1}{k} \sum_{m=1}^k F( e^{2\pi i m/k })
=
\int_{{S^1}} F(\varphi) \, d\varphi
+
O\biggl(\frac{1}{k}\biggr).
\label{eqdistdeltaktointwithkok}
\end{equation}
We now combine \eqref{eqtrFexplicit} and
\eqref{eqdistdeltaktointwithkok} and get
%
\begin{eqnarray}\label{eqasymptoticfortrace}
\tr(F)(\sigma)
&=&
\sum_{k=1}^N C_k \biggl(k \int_{{S^1}} F(\varphi) \, d\varphi+O(1) \biggr)
\nonumber
\\[-8pt]
\\[-8pt]
\nonumber
&=&
N \int_{{S^1}} F(\varphi) \, d\varphi+ O(T(\sigma)),
\end{eqnarray}
where we have used that $\sum_{k=1}^N kC_k = N$ and
$\sum_{k=1}^N C_k = T(\sigma)$. Notice that
\eqref{eqasymptoticfortrace} is independent of any
probability measure on $\SN$. Using the binomial theorem and
the fact that $0<T(\sigma)/N \leq1$ for all $\sigma$, we get
%
\begin{equation}
\frac{1}{N^d} (\tr(F)(\sigma))^d
=
\biggl(\int_{{S^1}} F(\varphi) \, d\varphi\biggr)^d
+ O_{F,d} \biggl( \frac{T(\sigma)}{N} \biggr),
\end{equation}
where the constant implicit in the big-O is independent of
$\sigma$ and $N$. We apply $\mathbb{E}_{\Theta}[\cdot]$ on both sides,
and this
proves the first part of the theorem.

The last statement follows from~\cite{NikZei}, Theorem 4.2, where it
is shown that if $g_\Theta(t)
\in\mathcal{F}(\vartheta,r,K)$ or $g_\Theta(t) \in
e\mathcal{F}(\vartheta,r,\gamma)$ then $\mathbb{E}_{\Theta
}[T(\sigma)] \sim
\vartheta\log(N)$.
\end{pf}
%

\begin{remark*}
In fact, for many probability distributions on $\SN$,
$\mathbb{E}_{\Theta}[T(\sigma)] = o(N)$. The only way for this not
to be true
is for $\sigma$ to frequently have only small cycles, which
will occur if $\Theta= (\theta_k)_{k=1}^\infty$ is
a sequence tending to zero very rapidly.
\end{remark*}

\section{Wreath product, traces and the generalized Ewens measure}
\label{secgeneralcase}
In this section we consider the traces of the wreath product
$S^1 \wr\SN$ (see, e.g.,~\cite{Wie03}). More precisely,
we consider random matrices of the form
\[
M(\sigma,z_1,\ldots,z_N)
:= \operatorname{diag}(z_1,\ldots,z_N)\cdot\sigma,
\]
where $\sigma$ is a random permutation of $\SN$, and $(z_j)_{j
\geq1}$ is a sequence of i.i.d. random variables with values
in $S^1$ (the complex unit circle), independent of $\sigma$.
Many groups closely related to $\SN$ give similar matrices, for
instance, the Weyl group of $\operatorname{SO}(2N)$.

The trace of a function $F$ is then extended in the obvious way by
%
\begin{equation}
\tr(F) = \tr(F,z_1,\ldots,z_N)(\sigma):= \sum_{k=1}^N F(\omega_k),
\end{equation}
where $(\omega_k)_{k=1}^N$ are the $N$ eigenvalues of
$M(\sigma,z_1,\ldots,z_N)$.

We now give a more explicit expression of $\tr(F)$.
%
\begin{lemma}\label{lembahroumi}
%
\begin{equation}
\label{eqtraceFwreath}
\tr(F) \stackrel{d}{=} \sum_{k=1}^N \sum_{m=1}^{C_k} k\Delta_k(F,Z_{k,m}),
\end{equation}
where $(Z_{k,m})_{k, m \geq1}$ is a sequence of independent random
variables which is independent of $(C_k)_{k \geq1}$ (the sequence of
cycle numbers of $\sigma$), with $Z_{k,m}$
equal in distribution to $\prod_{j=1}^k z_j$, and
%
\begin{equation}
\Delta_k(F,y) := \frac{1}{k} \sum_{\omega^k = y} F(\omega).
\end{equation}
\end{lemma}

\begin{pf}
The characteristic polynomial of $M(\sigma,z_1,\ldots,z_N)$ with
$\sigma\in\SN$ with cycle type $\la$, is given by
%
\begin{equation}
\label{eqznwreath}
\det\bigl(1-x M(\sigma,z_1,\ldots,z_N)\bigr)
=
\prod_{k=1}^{N} \prod_{m=1}^{C_k}\Biggl(1-x^{k} \prod_{j=1}^{k}
z_{j}^{k,m} \Biggr),
\end{equation}
where the sequence $(z_{j}^{k,m})_{k,m,j}$ is the same sequence
as $(z_j)_{j=1}^N$, but with a different numeration and
ordering. [Note that this is why it is crucial that the $(z_j)$
are i.i.d.] The proof of \eqref{eqznwreath} is similar to the
proof of \eqref{eqdefznla} and we thus omit the details. The
lemma now follows immediately from \eqref{eqznwreath}.
\end{pf}

As in Section~\ref{sectraces}, we can compute the generating function
of $\tr(F)$.

\begin{lemma}\label{le5.2}
We define
%
\begin{equation}
\chi_k(s):= \mathbb{E}\bigl[e^{is k \Delta_k(F,Z_{k,m})}\bigr].
\end{equation}

We then have
%
\begin{equation}
\sum_{N=0}^\infty h_N\mathbb{E}_{\Theta}\bigl[e^{i s\tr(F)}\bigr] t^N
=
\exp\Biggl( \sum_{k=1}^\infty\frac{\theta_k}{k} \chi_k(s) t^k \Biggr).
\label{eqgeneratingtrwreath}
\end{equation}
\end{lemma}

\begin{remark*}
Note that $\chi_k(s)$ is independent of $m$ since $Z_{k,1}
\stackrel{d}{=} Z_{k,2} \stackrel{d}{=} \cdots\stackrel{d}{=}
Z_{k,m}$.
\end{remark*}

\begin{pf*}{Proof of Lemma~\ref{le5.2}}
We compute $\mathbb{E}[\exp(i s\tr(F))]$. For this we use the
independence of~$C_k$ and $\Delta_k$ to obtain
%
\begin{eqnarray}
\mathbb{E}\bigl[e^{i s\tr(F)}\bigr]
&=&
\mathbb{E}\Biggl[\prod_{k=1}^N \prod_{m=1}^{C_k}e^{is k \Delta_k(F,Z_{k,m})}\Biggr]
=
\mathbb{E}\Biggl[\prod_{k=1}^N \prod_{m=1}^{C_k} \chi_k(s) \Biggr]
\nonumber
\\[-8pt]
\\[-8pt]
\nonumber
&=&
\mathbb{E}\Biggl[\prod_{k=1}^N (\chi_k(s))^{C_k} \Biggr].
\end{eqnarray}
The theorem now follows immediately from Lemma~\ref{lemcycleindextheorem}.
\end{pf*}

\begin{definition}
Let
%
\begin{equation}
g_{\tr(F)}(t)
:=
\sum_{k=1}^\infty\frac{\theta_k}{k} \chi_k(s) t^k.
\end{equation}
\end{definition}

\begin{theorem} \label{thmfouriertrF}
Assume $\mathbb{E}[|\Delta_k(F,Z_{k,1})|] = O(k^{-1-\delta})$ for
some $0<\delta\leq1$,
and assume $g_\Theta(t)$ is in $e\mathcal{F}(r,\vartheta,\gamma)$,
where $e\mathcal{F}(r,\vartheta,\gamma)$
is given in Definition~$\ref{deffunctionclasseFalphargamma}$.
Then $g_{\tr(F)}(t)$ is in $e\mathcal{F}(r,\vartheta,\min{\{\gamma
,\delta\}})$ and
%
\begin{equation}\label{eqcharfuncTrFwreath}\qquad
\mathbb{E}_{\Theta}\bigl[e^{is\tr(F)}\bigr]
=
\exp\Biggl(\sum_{k=1}^\infty\frac{\theta_k}{k} \bigl(\chi_k(s)-1\bigr) r^k \Biggr)
+ O\bigl(N^{-\min{\{\gamma,\delta\}}}\log(N)\bigr)
\end{equation}
and as $N$ tends to infinity, $\tr(F)$ converges in law to the
random variable
%
\begin{equation}\label{eqTraceLimitY}
Y := \sum_{k = 1}^\infty\sum_{m=1}^{P_k} k \Delta_k(F,Z_{k,m}),
\end{equation}
where $(P_k)_{k \geq1}$ is a sequence of independent Poisson
random variables, independent of $(Z_{k,m})_{k, m \geq1}$, and such
that $P_k$ has parameter
$\theta_k r^k/k$. Here, the series defining $Y$ is a.s. absolutely convergent.
\end{theorem}

\begin{remark*}
By linearity of trace, if $F$ is Riemann integrable one can
always subtract a suitable constant to make $\int_{S^1}
F(\varphi) \,d\varphi= 0$, which ensures $\Delta_k(F,z) \to
0$ as $k\to\infty$.
\end{remark*}

\begin{remark*}
One should compare equation~\eqref{eqTraceLimitY} with
equation~\eqref{eqtraceFwreath}. The replacement of the
cycle counts $C_k$ with $P_k$ is indicative of Feller coupling
for the generalized Ewens measure.
\end{remark*}

\begin{pf*}{Proof of Theorem~\ref{thmfouriertrF}}
We have
%
\begin{equation}
g_{\tr(F)}(t)
=
\sum_{k=1}^\infty\frac{\theta_k}{k} \chi_k(s) t^k
=
g_\Theta(t) + \sum_{k=1}^\infty\frac{\theta_k}{k} \bigl(\chi_k(s)-1\bigr) t^k.
\end{equation}
We now have
%
\begin{equation}
\qquad|\chi_k(s)-1|
\leq
\mathbb{E}\bigl[\bigl|e^{is k \Delta_k(F,Z_{k,1})}-1\bigr|\bigr]
\leq
\mathbb{E}[( k |\Delta_k(F,Z_{k,1})| )]
=
O(k^{-\delta}).
\end{equation}
On the other hand, we have $\theta_k = O(r^{-k})$. This follows
immediately from the
fact that $g_\Theta(t)$ is in $e\mathcal{F}(r,\vartheta,\gamma)$.
We thus have
%
\begin{equation}
\frac{\theta_k}{k} \bigl(\chi_k(s)-1\bigr) = O(r^{-k} k^{-1-\delta}).
\end{equation}
This shows that $g_{\tr(F)}(t) \in e\mathcal{F}(r,\vartheta,\min\{
\gamma,\delta\})$.

Since $g_\Theta(t)\in e\mathcal{F}(r,\vartheta,\gamma)$, we can
write $g_\Theta(t) = \vartheta\log( \frac{1}{1-t/r} ) +
g_0(t)$ with $g_0(r) < \infty$. Thus,
%
\begin{equation}
g_{\tr(F)}(t) = \vartheta\log\biggl( \frac{1}{1-t/r} \biggr) + g_0(t) + \sum
_{k=1}^\infty\frac{\theta_k}{k} \bigl(\chi_k(s)-1\bigr) t^k.
\end{equation}
We get with Theorem~\ref{thmhwang2} that
%
\begin{equation}
\qquad h_N \mathbb{E}_{\Theta}\bigl[e^{is\tr(F)}\bigr]
=
\frac{N^{\vartheta-1}}{r^N\Gamma(\vartheta)} \exp\Biggl( g_0(r) + \sum
_{k=1}^\infty\frac{\theta_k}{k}\bigl (\chi_k(s)-1\bigr) r^k \Biggr) +R_N
\end{equation}
with
%
\begin{equation}
R_N =O\biggl(\frac{N^{\vartheta-1-\min\{\gamma,\delta\}}\log(N)}{r^N} \biggr).
\end{equation}
Dividing by $h_N$ proves equation~\eqref{eqcharfuncTrFwreath}.

Using the characteristic function of $\tr(F)$, we can deduce
its convergence in law to $Y$. The absolute convergence of the
series in \eqref{eqTraceLimitY} comes from
%
\begin{eqnarray}
\mathbb{E}\Biggl[\sum_{k=1}^\infty\sum_{m=1}^{P_k} k |\Delta
_k(F,Z_{k,m})| \Biggr]
& =& \sum_{k=1}^\infty k \mathbb{E}[P_k] \mathbb{E}[|\Delta
_k(F,Z_{k,1})|] \nonumber\\
& =& \sum_{k \geq1} \theta_k r^k O(k^{-1-\delta}) \\
& <& \infty, \nonumber
\end{eqnarray}
since $ \theta_k = O(r^{-k})$. Now, for $s \in\mathbb{R}$ and $k
\geq1$,
%
\begin{eqnarray}
\mathbb{E}\bigl[e^{is k \sum_{m=1}^{P_k} \Delta_k(F,Z_{k,m})} \bigr]
& =& \mathbb{E}\bigl[\bigl( \mathbb{E}\bigl[e^{is k \Delta_k(F,Z_{k,1})} \bigr] \bigr)^{P_k} \bigr]
= \mathbb{E}[(\chi_k(s) )^{P_k} ]
\nonumber
\\[-8pt]
\\[-8pt]
\nonumber
& =&
\exp\biggl( \frac{\theta_k}{k} \bigl(\chi_k(s)-1\bigr) r^k \biggr).
\end{eqnarray}
Hence, by absolute convergence,
%
\begin{equation}
\mathbb{E}[e^{is Y} ] = \exp\Biggl( \sum_{k=1}^\infty\frac{\theta_k}{k}
\bigl(\chi_k(s)-1\bigr) r^k \Biggr),
\end{equation}
and thus by equation~\eqref{eqcharfuncTrFwreath}, $
\mathbb{E}_{\Theta}[e^{is\tr(F)}] \to\mathbb{E}[e^{i s Y}]$ as
$N\to\infty$, and thus
$\tr(F)$ converges in law to $Y$.
\end{pf*}

With a more direct approach one can prove the convergence in
law of $\tr(F)$ to~$Y$ (albeit without a rate of convergence)
under slightly weaker conditions.
%
\begin{theorem} \label{poisson2}
Assume that $g_\Theta(t)$ is in
$e\mathcal{F}(r,\vartheta,\gamma)$ and that
%
\begin{equation}
\sum_{k=1}^\infty k^{(1- \vartheta)_+} \mathbb{E}[|\Delta_k(F,Z_{k,1})|]
< \infty.
\end{equation}
Then $\tr(F)$ converges in law to $Y$, where $Y$ is given by
\eqref{eqTraceLimitY}.
\end{theorem}

\begin{pf}
Under these conditions, the absolute convergence of the series
defining $Y$ is checked as follows:
%
\begin{eqnarray}
\mathbb{E}\Biggl[\sum_{k=1}^\infty\sum_{m=1}^{P_k} k |\Delta
_k(F,Z_{k,m})| \Biggr]
& =& \sum_{k=1}^\infty k \mathbb{E}[P_k] \mathbb{E}[|\Delta
_k(F,Z_{k,1})|] \\
&= &\sum_{k=1}^\infty\theta_k r^k \mathbb{E}[|\Delta_k(F,Z_{k,1})|]
\\
&=& O \Biggl( \sum_{k=1}^\infty\mathbb{E}[|\Delta_k(F,Z_{k,1})|] \Biggr)
\end{eqnarray}
and this converges by assumption.

In~\cite{NikZei}, Corollary 3.1.1, it is proven that for all
fixed $b \geq1$, $(C_1,C_2, \ldots, C_b)$ tends in law to
$(P_1,P_2, \ldots, P_b)$ when the dimension $N$ goes to
infinity.

Now let
%
\begin{equation}
\tr_b(F) := \sum_{k=1}^b \sum_{m=1}^{C_k} k \Delta_k(F,Z_{k,m})
\end{equation}
and
%
\begin{equation}
Y_b := \sum_{k=1}^b \sum_{m=1}^{P_k} k \Delta_k(F,Z_{k,m}).
\end{equation}
The same argumentation as in Theorem~\ref{thmfouriertrF} gives
%
\begin{equation}
\sum_{N=0}^\infty h_N \mathbb{E}_{\Theta}\bigl[e^{i s \tr_b(F) } \bigr] t^N =
\exp\Biggl( \sum_{k=1}^b \frac{\theta_k}{k} \bigl( \chi_k(s) - 1\bigr) t^k \Biggr)
e^{g_\Theta(t)}
\end{equation}
from which follows (again by the same reasoning as in
Theorem~\ref{thmfouriertrF})
%
\begin{equation}
\mathbb{E}_{\Theta}\bigl[e^{i s \tr_b(F) } \bigr] \to\mathbb{E}[e^{i s Y_b
} ]
\end{equation}
as $N\to\infty$. Here $b$ is fixed but arbitrary, so this
convergence implies (by using the inequality $|e^{ix} - e^{iy}|
\leq|x-y|$) that
%
\begin{eqnarray}
&&\limsup_{N\to\infty} \bigl| \mathbb{E}_{\Theta}\bigl[e^{ i s \tr(F)}\bigr] -
\mathbb{E}[e^{ i s Y}] \bigr| \nonumber\\
&&\qquad\leq|s| \limsup_{N\to\infty}
\sum_{k=b+1}^{\infty} k \mathbb{E}[|\Delta_k(F,Z_{k,1}) |]
\mathbb{E}'[ (C_k + P_k)] \\
&&\qquad\leq|s| \sum_{k=b+1}^{\infty} k \mathbb{E}
[|\Delta_k(F,Z_{k,1}) |] H_k,\nonumber
\end{eqnarray}
where $\mathbb{E}'$ is the expectation over the product measure
of $\mathbb{P}_\Theta$ and the measures occurring from
$(P_k)_{k\geq1}$, and where
%
\begin{equation}
H_k = \mathbb{E}[P_k] +
\sup_{N \geq1} \mathbb{E}_{\Theta}[C_{k}].
\end{equation}
Therefore, the theorem is proven if we show that
%
\begin{equation}
\sum_{k=1}^\infty k H_k \mathbb{E}[|\Delta_k(F,Z_{k,1}) |] <
\infty.
\end{equation}

Ercolani and Ueltschi~\cite{cyclestructureueltschi}, Proposition
2.1(c), show that
%
\begin{equation}
\mathbb{E}_{\Theta}[C_{k}] =
\cases{
\displaystyle\frac{\theta_k}{k} \frac{h_{N-k}}{h_N}, & \quad$\mbox{ if } k\leq
N$,
\vspace*{2pt}\cr
0, &\quad $\mbox{ if } k>N.$
}
\end{equation}

By Lemma~\ref{lemasymptotichn2}, we have, for some $A > 0$
and for $N$ going to infinity,
%
\begin{equation}
h_N \sim A (N+1)^{\vartheta-1}/ r^N,
\end{equation}
so
%
\begin{equation}
\frac{h_{N-k}}{h_N} = O\biggl(r^{k} \biggl(1-\frac{k}{N+1}\biggr)^{\vartheta-1}\biggr).
\end{equation}
Now, for $k$ fixed,
%
\begin{eqnarray}
&&\max_{N\geq k} \biggl(1-\frac{k}{N+1}\biggr)^{\vartheta-1}
\nonumber
\\[-8pt]
\\[-8pt]
\nonumber
&&\qquad =
\cases{
1, &\quad$\mbox{if } \vartheta\geq1$,\vspace*{2pt}\cr
(k+1)^{1-\vartheta}, & \quad$\mbox{if $\vartheta< 1$ (attained at
$N =
k$)}$.}
\end{eqnarray}
Since $g_\Theta(t) \in e\mathcal{F}(r,\vartheta,\gamma)$, we
have $\theta_k r^k = O(1)$, and so we deduce that
%
\begin{equation}
\sup_{N \geq1} \mathbb{E}_{\Theta}[C_k] = O\biggl(\frac{\theta_k}{k}
r^k (k+1)^{(1-\vartheta)_+} \biggr) = O\bigl(k^{-1+(1-\vartheta)_+ }\bigr).
\end{equation}

Finally, since
%
\begin{equation}
\mathbb{E}[P_k] = \theta_k r^k / k = O(1/k)
\end{equation}
we have
%
\begin{equation}
H_k = O\bigl(k^{-1+(1-\vartheta)_+ }\bigr),\vadjust{\goodbreak}
\end{equation}
and so
%
\begin{equation}
\qquad\sum_{k=1}^\infty k H_k \mathbb{E}[|\Delta_k(F,Z_{k,1}) |]
= O \Biggl( \sum_{k=1}^\infty k^{(1-\vartheta)_+} \mathbb{E}[|\Delta
_k(F,Z_{k,1}) |] \Biggr) < \infty
\end{equation}
as required.
\end{pf}

\begin{remark*}
In the case when $z_i$ are all equal to $1$ almost surely, then $\Delta
_k(F,Z_{k,m}) = \Delta_k(F)$
as given in \eqref{eqDeltak}, and we are back in the case of
permutation matrices. Thus these two theorems fulfill the
promise made in Section~\ref{sectraces}.
\end{remark*}

The following result gives sufficient conditions, expressed
only in terms of the Fourier coefficients of $F$, under which
we can be assured the conditions of
Theorems~\ref{thmfouriertrF} and~\ref{poisson2} apply.
%
\begin{theorem} \label{fourier}
Let us suppose that $F$ is continuous and for $m \in
\mathbb{Z}$, let us define the $m$th Fourier coefficient of
$F$ by
%
\begin{equation}
c_m(F) := \frac{1}{2\pi} \int_0^{2 \pi}
e^{-i m x} F(e^{i x}) \,d x.
\end{equation}
We assume that the mean value of $F$ vanishes, that is, $c_0(F) = 0$.
If for some $\delta\in(0,1]$, $c_m (F) = O(|m|^{-1-\delta})$
when $|m|$ goes to infinity then
%
\begin{equation}
\mathbb{E}[|\Delta_k(F,Z_{k,1})|] = O(k^{-1-\delta}).
\end{equation}

If there exists $s >(1-\vartheta)_+$ such that
%
\begin{equation}
\sum_{m \in\Z} |m|^{s} |c_m(F)| < \infty
\end{equation}
then
%
\begin{equation}
\sum_{k=1}^\infty k^{(1- \vartheta)_+} \mathbb{E}[|\Delta
_k(F,Z_{k,1})|] < \infty.
\end{equation}
\end{theorem}

\begin{remark*}
If the assumptions of Theorem $\ref{fourier}$ are satisfied, except
that $c_0(F) = 0$,
then one can still apply the result to the function $F - c_0(F)$, and
deduce, from Theorem
$\ref{thmfouriertrF}$ or Theorem $\ref{poisson2}$, that $\tr(F) - N
c_0(F)$ converges in law to $Y$, where $Y$ is given by \eqref{eqTraceLimitY}.
\end{remark*}

\begin{pf*}{Proof of Theorem~\ref{fourier}}
Since $F$ is continuous, one has, for all $x \in[0, 2 \pi)$,
%
\begin{equation}
F(e^{i x}) = \lim_{n \rightarrow\infty} \sum_{m \in
\mathbb{Z}} \frac{(n-|m|)_+}{n} c_m(F) e^{i m x},
\end{equation}
by using the Fej\'er kernel. Now, by assumption,
%
\begin{equation}
\sum_{m \in\mathbb{Z}} |c_m(F)| < \infty,
\end{equation}
and hence, by dominated convergence,
%
\begin{equation}
F(e^{i x}) = \sum_{m \in\mathbb{Z}} c_m(F) e^{i m x},
\end{equation}
where the series is absolutely convergent. Since $c_0(F) = 0$,
one deduces that for all $k \geq1$ and $x \in[0, 2
\pi)$,
\begin{eqnarray*}
\Delta_{k} (F, e^{i x})
& =&
\frac{1}{k} \sum_{j=0}^{k-1} F \bigl(e^{i(x
+ 2 j \pi)/k} \bigr) \\
& =& \frac{1}{k} \sum_{m \in\Z\setminus\{0\}} c_m(F) \Biggl( \sum_{j=0}^{k-1}
e^{i m (x + 2 j \pi)/k} \Biggr) \\
& =& \mathop{\sum_{m \in\Z\setminus\{0\},}}_{ k | m} c_{m} (F) e^{i
m x/ k}.
\end{eqnarray*}

If $F$ satisfies the first assumption, that $c_m (F) =
O(|m|^{-1-\delta})$, then
%
\begin{eqnarray}
\sup_{x \in[0,2\pi)} |\Delta_{k} (F,e^{i x})|
&\leq&\mathop{\sum_{m \in\Z\setminus\{0\},}}_{ k | m} |c_{m} (F)|
= O \biggl( \mathop{\sum_{m \in\Z\setminus\{0\},}}_{ k | m}
|m|^{-1-\delta} \biggr) \\
&=& O(k^{-1-\delta})
\end{eqnarray}
for $k$ going to infinity, which clearly implies
$\mathbb{E}[|\Delta_k(F,Z_{k,1})|] = O(k^{-1-\delta})$.

If $F$ satisfies the second assumption, one has
%
\begin{eqnarray}
\qquad\sum_{k =1}^\infty k^{(1-\vartheta)_+} \sup_{x \in[0,2\pi)}
|\Delta_{k} (F,e^{ix})|
&\leq&
\sum_{k =1}^\infty
k^{(1-\vartheta)_+} \mathop{\sum_{m \in\Z\setminus\{0\}, }}_{k |
m} |c_{m} (F)|\\
&\leq&
\sum_{k =1}^\infty
\mathop{\sum_{m \in\Z\setminus\{0\}, }}_{k | m} |m|^{(1-\vartheta
)_+} |c_{m} (F)|\\
&\leq&
\sum_{m \in\Z\setminus\{0\}} |c_{m} (F)| |m|^{(1-\vartheta)_+}
\tau(|m|),
\end{eqnarray}
where $\tau(|m|)$ denotes the number of divisors of $|m|$. Since $\tau
(|m|) = O(|m|^{\varepsilon})$ for all $\varepsilon> 0$, one
deduces that
%
\begin{equation}
\sum_{k \geq1} k^{(1-\vartheta)_+} \mathbb{E}[|\Delta_{k} (F,Z_{k,1})|]
=
O \biggl(\sum_{m \in\Z\setminus\{0\}} |c_{m} (F)| |m|^{s} \biggr) < \infty.
\end{equation}
The proof of the theorem is complete.
\end{pf*}

%
\begin{corollary} \label{sobolev}
Let $F$ be a continuous function from $S^1$ to
$\C$, contained in a Sobolev space $H^{s}$ for some $s
> 1/2 + (1-\vartheta)_+$.
Then, the second condition of Theorem~$\ref{fourier}$ is
fulfilled, and thus also the conditions of
Theorem~$\ref{poisson2}$.
\end{corollary}
\begin{pf}
By the Cauchy--Schwarz inequality, one has, for any $\alpha\in
((1-\vartheta)_+, s-1/2)$,
%
\begin{eqnarray}
&&\sum_{m \in\Z\setminus\{0\}} |m|^{\alpha} |c_m(F)|
\nonumber
\\[-8pt]
\\[-8pt]
\nonumber
&&\qquad\leq
\biggl(\sum_{m \in
\Z\setminus\{0\}} |m|^{2s} |c_m(F)|^2 \biggr)^{1/2}
\biggl( \sum_{m \in\Z\setminus\{0\}} |m|^{2(\alpha-s)} \biggr)^{1/2},
\end{eqnarray}
which is finite since $F \in H^{s}$ and $2(\alpha-s) < -1$.
\end{pf}

\begin{remark*}
Note that it is not always obvious to estimate directly the
Fourier coefficients of a function $F$, however, standard
results from Fourier analysis concerning the differentiability
of $F$ yield sufficient bounds on the decay of the Fourier
coefficients of $F$ for the conditions of Theorem~$\ref{fourier}$
to be checked (see, e.g,~\cite{Kor}, Chapter 9).
\end{remark*}

\section{Diverging variance for the classical Ewens measure}\label
{sectEwensDivVar}

In the previous two sections, we have been considering the
convergence of $\tr(F)$ to some limit for random permutation
matrices (and their generalization to wreath products), where
the underlying probability space is the generalized Ewens
measure. The conditions we have used have all implied that the
variance of $\tr(F)$ stays bounded as $N\to\infty$.

A recent paper by Ben Arous and Dang~\cite{benarousdang}
dealing with $\tr(F)$ for real $F$ and for random permutation matrices
in the
special case of the classical Ewens measure (which is when
$\theta_k = \theta$, a constant), demonstrates a dichotomy
between converging and diverging variance for $\tr(F)$ in the
classical Ewens distribution. In the former case they also show
convergence of $\tr(F)$ to an explicit finite limit, and in the
latter case they prove the following central limit theorem.
%
\begin{theorem}[(Ben Arous and Dang)]\label{thmBenArous}
Let $F\dvtx\C\to\R$ be given and assume that
%
%
\begin{equation}
V_N := \theta\sum_{k=1}^N k \Delta_k(F)^2
\end{equation}
tends to infinity as $N\to\infty$
and
%
\begin{equation}\label{equnnecessarycondition}
\max_{1 \leq k \leq N} k |\Delta_k(F)| = o\bigl(\sqrt{V_N}\bigr)
\end{equation}
then,
%
\begin{equation}
\frac{\tr(F)-\mathbb{E}[\tr(F)]}{\sqrt{V_N}} \stackrel
{d}{\longrightarrow} \mathcal{N}(0,1) .
\end{equation}
\end{theorem}

In the generalized Ewens measure, we are currently unable to
apply the function theoretic methods to prove weak convergence
results in the case of diverging variance. However, for the
classical Ewens measure we are able to prove a similar central
limit theorem the wreath product, with slightly extended
application, in the sense that condition
\eqref{equnnecessarycondition} can be weakened from a
sup-norm to a $p$-norm.

\begin{theorem}
\label{thmtraceinfinitecase}
Let $F\dvtx\C\to\R$ be given and assume that
%
\begin{equation}\label{eqVN}
V_N := \theta\sum_{k=1}^N k \mathbb{E}[(\Delta_k(F,Z_{k,1}))^2]
\end{equation}
tends to infinity as $N\to\infty$.
Assume further that there exists a $p>\max\{ \frac{1}{\theta} , 2 \}
$ such
that
%
\begin{equation}\label{eqinfvarcond}
\sum_{k=1}^{N} k^{p-1} \mathbb{E}[|\Delta_k(F,Z_{k,1})|^p] =
o(V_N^{p/2})
\end{equation}
with $\Delta_k(F,z) = \frac{1}{k} \sum_{\omega^k = z} F(\omega)$.
Then
%
\begin{equation}
\biggl( \frac{\tr(F) - E_N}{ \sqrt{V_N}} \biggr)_{N \geq1}
\end{equation}
converges in distribution to a standard Gaussian random
variable, where
%
\begin{equation}\label{eqEN}
E_N := \theta\sum_{k = 1}^N \mathbb{E}[\Delta_k(F,Z_{k,1})].
\end{equation}
\end{theorem}

The behavior for complex functions $F$ can be computed in a similar way.
We consider here only real $F$ to keep the notation simple and to avoid
further technicalities.

\begin{remark*}
Recall that without loss of generality we may assume $F$ has
mean zero in the sense that $\int_{0}^{2\pi} F(e^{i x}) \,dx =
0$. We remark that\vadjust{\goodbreak} this does not necessarily imply that
$\Delta_k(F,z)$ tends to zero, even though $\Delta_k(F,z)$ is a
discretization of the integral, without the assumption of
further smoothness conditions. Moreover, note that in the framework of
the symmetric
group (i.e., $Z_{k,1}= 1$ almost surely), the assumption
$(\ref{eqinfvarcond})$ is implied by the condition $(\ref
{equnnecessarycondition})$ given
in~\cite{benarousdang}. Indeed, if (\ref{equnnecessarycondition}) is
satisfied, then for
$p >2$, one has
\begin{eqnarray*}
\sum_{k=1}^N k^{p-1} |\Delta_k(F)|^p & \leq&
\Bigl( \max_{1 \leq k \leq N} (k |\Delta_k(F)|) \Bigr)^{p-2}
\sum_{k=1}^N k |\Delta_k(F)|^2 \\
& =& o\bigl(V_N^{(p-2)/2}\bigr) O(V_N) = o(V_N^{p/2}).
\end{eqnarray*}
\end{remark*}

In the proof of Theorem~\ref{thmtraceinfinitecase}, we use the Feller
coupling, which allows the random
variables $C_k$ and $P_k$ to be defined on the same space
and to replace the weak convergence $C_k \stackrel{d}{\to} P_k$ by
convergence in probability (but not a.s. convergence).
This coupling exists only for the classical Ewens measure and
thus $P_k$ are independent Poisson distributed random variables with
$\mathbb{E}[P_k] = \frac{\theta}{k}$.
The construction and further details can be found, for instance, in
\cite{ABT02}, Sections 1 and 4.

The Feller coupling allows us to prove Theorem~\ref
{thmtraceinfinitecase} with $C_k$ replaced by
$P_k$, and the following lemma allows us to estimate the
distance between the two.
%
\begin{lemma}[(Ben Arous and Dang~\cite{benarousdang})]
\label{lemfellerBenDang} For any $\theta>0$ there exists a
constant $K(\theta)$ depending on $\theta$, such that for every
$1\leq m\leq N$,
%
\begin{equation}
\mathbb{E}[|C_k -P_k |]\leq\frac{K(\theta)}{N}+\frac\theta N\Psi_N(k),
\end{equation}
where
%
\begin{equation}
\Psi_N(k):=\pmatrix{N-k+\theta-1\cr N-k}\pmatrix{N+\theta-1\cr N}^{-1}.
\end{equation}
\end{lemma}

\begin{pf*}{Proof of Theorem~\ref{thmtraceinfinitecase}}
The main idea of the proof is to define the auxiliary random
variable
%
\begin{equation}
Y_N(F)
:=
\sum_{k=1}^N \sum_{m=1}^{P_k} k\Delta_k (F,Z_{k,m})
\end{equation}
and to show that $\tr(F)$ and $Y_N(F)$ have the same asymptotic
behavior after normalization, and that (again after
normalization) $Y_N(F)$ satisfies a central limit theorem.

First, we will show that
%
\begin{equation}
\mathbb{E}[|\tr(F) -Y_N(F)|] = o((V_N)^{{1}/{2}}).
\label{eqtr-trH=osigma}
\end{equation}

We use Lemma~\ref{lemfellerBenDang} and that $Z_{k,m}$ are
independent of $C_k$ and $P_k$ to get
%
\begin{eqnarray} \label{eqboundingTrFYF}
&&\mathbb{E}[|\tr(F) - Y_N(F)|]\nonumber\\
&&\qquad\leq
\mathbb{E}\Biggl[\sum_{k=1}^N |C_k-P_k| k \mathbb{E}[|\Delta_k
(F,Z_{k,1})|]\Biggr]\\
&&\qquad\leq
\frac{K(\theta)}{N} \sum_{k=1}^N k\mathbb{E}[|\Delta
_k(F,Z_{k,1})|] +\frac{\theta}{N} \sum_{k=1}^N k\mathbb{E}[|\Delta
_k(F,Z_{k,1})|] \Psi_N(k).\nonumber
\end{eqnarray}

For the first term, we apply Jensen's inequality and condition
\eqref{eqinfvarcond} to obtain
%
\begin{eqnarray}
\frac1N \sum_{k=1}^N k\mathbb{E}[|\Delta_k(F,Z_{k,1})|] &\leq&\Biggl(
\frac{1}{N} \sum_{k=1}^N k^p \mathbb{E}[|\Delta_k(F,Z_{k,1})|]^p
\Biggr)^{1/p} \\
&\leq&\Biggl( \sum_{k=1}^N k^{p-1} \mathbb{E}[|\Delta_k(F,Z_{k,1})|^p]
\Biggr)^{1/p} \\
&= &o(V_N^{1/2}) .
\end{eqnarray}

Now we deal with the second term in \eqref{eqboundingTrFYF}.
If $\theta\geq1$ then $\Psi_N(k)$ is bounded by~$1$, so the
same argument as above shows that the second summand is also
$o(V_N^{1/2})$ in this case.

If $0<\theta< 1$, we have to be more careful. A simple
computation shows that there exists constants $K_1,K_2$ such
that
%
\begin{equation}\label{equpperboundbinom}
\Psi_N(k)
\leq
\cases{
\displaystyle K_1 \biggl(1-\frac{k}{N}\biggr)^{\theta-1}, &\quad $\mbox{for }k< N,$
\vspace*{2pt}\cr
\displaystyle K_2 N^{1-\theta} , &\quad $\mbox{for } k=N$}
\end{equation}
and so
%
\begin{eqnarray}
&&\frac{\theta}{N} \sum_{k=1}^N k\mathbb{E}[|\Delta_k(F,Z_{k,1})|]
\Psi_N(k)
\\
&&\qquad\leq
\theta K_2 N^{1-\theta} \mathbb{E}[|\Delta_N(F,Z_{N,1})|]
\nonumber
\\[-8pt]
\\[-8pt]
\nonumber
&&\qquad\quad{}+
\frac{\theta K_1}{N} \sum_{k=1}^{N-1} k\mathbb{E}[|\Delta
_k(F,Z_{k,1})|] \biggl(1-\frac{k}{N}\biggr)^{\theta-1}.
\end{eqnarray}

Using the value of $p$ given in the conditions of the theorem,
%
\begin{eqnarray}
N^{1-\theta} \mathbb{E}[|\Delta_N(F,Z_{N,1})|] &= &( N^{p-p\theta}
\mathbb{E}[|\Delta_N(F,Z_{N,1})|]^p )^{1/p} \\
&\leq&\Biggl( \sum_{k=1}^N k^{p-p\theta} \mathbb{E}[|\Delta
_k(F,Z_{k,1})|]^p \Biggr)^{1/p} \\
&\leq&\Biggl( \sum_{k=1}^N k^{p-1} \mathbb{E}[|\Delta_k(F,Z_{k,1})|^p] \Biggr)^{1/p}
\end{eqnarray}
since $k^{p-p\theta} \leq k^{p-1}$ (since $p\theta>1$) and
$\mathbb{E}[|\Delta_k(F,Z_{k,1})|]^p \leq\mathbb{E}[|\Delta
_k(F,Z_{k,1})|^p]$
(since $p>1$). Thus, by condition~\eqref{eqinfvarcond}, this
is $o(V_N^{1/2})$.

H\"older's inequality gives
%
\begin{eqnarray}
\qquad&&\frac{1}{N} \sum_{k=1}^{N-1} k\mathbb{E}[|\Delta_k(F,Z_{k,1})|]
\biggl(1-\frac{k}{N}\biggr)^{\theta-1}
\nonumber
\\[-8pt]
\\[-8pt]
\nonumber
&&\qquad\leq\Biggl( \frac{1}{N} \sum_{k=1}^{N-1}
k^p\mathbb{E}[|\Delta_k(F,Z_{k,1})|^p] \Biggr)^{1/p} \Biggl( \frac{1}{N}
\sum_{k=1}^{N-1} \biggl(1-\frac{k}{N}\biggr)^{q(\theta-1)}
\Biggr)^{1/q}
\end{eqnarray}
with $\frac{1}{p} +\frac{1}{q}=1$.

After re-ordering the sum, the second factor is
$(N^{q(1-\theta)-1} \sum_{j=1}^{N-1} j^{-q(1-\theta)})^{1/q}$, and if
$q(\theta-1)> -1$ then it is bounded above by a constant. Note
that
\begin{eqnarray*}
(\theta-1)q > -1
\quad&\Longleftrightarrow&\quad
(1-\theta) < \frac{1}{q}
\quad\Longleftrightarrow\quad
(1-\theta) < 1- \frac{1}{p}\\
&\Longleftrightarrow&\quad
\theta> \frac{1}{p}
\end{eqnarray*}
and thus condition~\eqref{eqinfvarcond} now ensures the
existence of a $p> \frac{1}{\theta}$ such that the first factor
is $o((V_N)^{{1}/{2}})$ and the second factor is
bounded. This proves \eqref{eqtr-trH=osigma}.

Therefore Slutsky's theorem implies that $Y_N(F)$ and $\tr(F)$
have the same asymptotic distribution after scaling. Thus it
suffices to show that
%
\begin{equation}
\frac{Y_N(F)- E_N}{\sqrt{V_n}}
\end{equation}
converges in law to a standard Gaussian random variable.

We calculate the mean of $Y_N(F)$ by first taking expectation
with respect to $Z_{k,m}$ and then with respect to $P_k$, to
obtain
%
\begin{eqnarray}
\mathbb{E}[Y_N(F)]& =& \sum_{k=1}^N \mathbb{E}\Biggl[\sum_{m=1}^{P_k} k
\mathbb{E}[\Delta_{k}(F,Z_{k,m})]\Biggr]
\nonumber
\\[-8pt]
\\[-8pt]
\nonumber
&=& \sum_{k=1}^N \mathbb{E}[P_k] k
\mathbb{E}[\Delta_{k}(F,Z_{k,1})],
\end{eqnarray}
where we use the fact that $\mathbb{E}[\Delta_{k}(F,Z_{k,m})] =
\mathbb{E}[\Delta_{k}(F,Z_{k,1})]$ for all $m$. Finally, since
$\mathbb{E}[P_k]
= \theta/k$, we see that $\mathbb{E}[Y_N(F)] = E_N$ as defined in
\eqref{eqEN}.

For the variance, since $P_k$ and $Z_{k,m}$ are all
independent, one can move the sum outside the variance, to
obtain
%
\begin{equation}
\var(Y_N(F)) = \sum_{k=1}^N \var\Biggl( \sum_{m=1}^{P_k} k \Delta
_{k}(F,Z_{k,m}) \Biggr).
\end{equation}
Now, the variance of a sum of random length of i.i.d. random
variables is given by the following formula:
%
\begin{equation}
\var\Biggl( \sum_{m=1}^P X_m \Biggr) = \var(X_1) \mathbb{E}[P] + \var(P)
\mathbb{E}[X_1]^2,
\end{equation}
if $(X_m)_{m \geq1}$ are i.i.d., $L^2$ random variables, and if $P$ is
an $L^2$ variable,
independent of $(X_m)_{m \geq1}$
(this result can be proved by a straightforward calculation). Letting
$X_m = k \Delta_{k}(F,Z_{k,m})$ and $P = P_k$ and knowing that
$\mathbb{E}[P_k] = \var(P_k) = \theta/k$, we deduce that $\var
(Y_N(F))= V_N$.

Finally we apply the Lyapunov central limit theorem since
$Y_N(F)$ is a sum of independent random variables. We will show
that
%
\begin{eqnarray}\label{eqLyapunov}
&&\sum_{k=1}^N \mathbb{E}\Biggl[\Biggl|\sum_{m=1}^{P_k} k \Delta_k(F,Z_{k,m})
-\mathbb{E}[k P_k \Delta_k (F,Z_{k,1})] \Biggr|^p\Biggr]
\nonumber
\\[-8pt]
\\[-8pt]
\nonumber
&&\qquad\ll\sum_{k=1}^N
k^{p-1} \mathbb{E}[| \Delta_k (F,Z_{k,1}) |^p ]
\end{eqnarray}
and by condition~\eqref{eqinfvarcond}, with $p>2$, this is
$o(V_N^{p/2})$ which means $\frac{Y_N(F)- E_N}{\sqrt{V_n}}$
converges in law to a standard Gaussian random variable.

For simplicity, let $P$ be a Poisson random variable with
parameter $\theta/k$ (we think of $k$ as being large), and let
$X_m = k \Delta_k (F,Z_{k,1})$ be i.i.d. random variables with
$\mathbb{E}[|X_m|^p]$ finite. Then
%
\begin{eqnarray}\label{eqlyapunovcenteredsums}
&&\mathbb{E}\Biggl[\Biggl|\sum_{m=1}^P X_m - \mathbb{E}[P] \mathbb{E}[X_1] \Biggr|^p\Biggr]\nonumber\\
&&\qquad= \mathbb{E}\Biggl[\Biggl|\sum_{m=1}^P (X_m - \mathbb{E}[X_1]) + (P-\mathbb
{E}[P]) \mathbb{E}[X_1] \Biggr|^p\Biggr] \\
&&\qquad\leq\Biggl(\mathbb{E}\Biggl[\Biggl|\sum_{m=1}^P (X_m - \mathbb{E}[X_1])\Biggr|^p\Biggr]^{1/p} +
\mathbb{E}\bigl[|P-\mathbb{E}[P]|^p\bigr]^{1/p} \mathbb{E}[X_1]
\Biggr)^p\nonumber
\end{eqnarray}
by the generalized triangle inequality.

Now, for all $p>1$,
%
\begin{equation}
\mathbb{E}\bigl[|P-\mathbb{E}[P]|^p\bigr] \leq\mathbb{E}[P^p] \ll\mathbb{E}[P]
\end{equation}
as $\mathbb{E}[P] \to0$, and so the second term in
\eqref{eqlyapunovcenteredsums} is
%
\begin{eqnarray}
\mathbb{E}\bigl[|P-\mathbb{E}[P]|^p\bigr]^{1/p} \mathbb{E}[X_1] &\ll&\mathbb
{E}[P]^{1/p} \mathbb{E}[X_1] \\
&\ll&( \mathbb{E}[P] \mathbb{E}[|X_1|^p] )^{1/p}
\end{eqnarray}
by H\"older's inequality, since $p>1$.

To bound the first term in \eqref{eqlyapunovcenteredsums},
let $q_n = \mathbb{P}[P = n]$, and note that
%
\begin{eqnarray}
\mathbb{E}\Biggl[\Biggl|\sum_{m=1}^P (X_m - \mathbb{E}[X_1])\Biggr|^p\Biggr] &=& \sum
_{n=0}^\infty q_n \mathbb{E}\Biggl[\Biggl| \sum_{m=1}^n (X_m - \mathbb{E}[X_1])
\Biggr|^p \Biggr] \\
&\leq&\sum_{n=0}^\infty q_n \Biggl( \sum_{m=1}^n \mathbb{E}\bigl[|X_m - \mathbb
{E}[X_1]|^p \bigr]^{1/p} \Biggr)^p \\
&= &\sum_{n=0}^\infty q_n n^p \mathbb{E}\bigl[|X_1 - \mathbb{E}[X_1]|^p \bigr]
\\
&=& \mathbb{E}[P^p] \mathbb{E}\bigl[|X_1 - \mathbb{E}[X_1]|^p \bigr] \\
&\ll&\mathbb{E}[P] \mathbb{E}[|X_1|^p ].
\end{eqnarray}

Thus,
%
\begin{eqnarray}
&&\mathbb{E}\Biggl[\Biggl|\sum_{m=1}^P X_m - \mathbb{E}[P] \mathbb{E}[X_1] \Biggr|^p\Biggr]
\nonumber\\
&&\qquad\ll\bigl( ( \mathbb{E}[P] \mathbb{E}[|X_1|^p ] )^{1/p} + ( \mathbb
{E}[P] \mathbb{E}[|X_1|^p] )^{1/p} \bigr)^p \\
&&\qquad\ll\mathbb{E}[P] \mathbb{E}[|X_1|^p ].\nonumber
\end{eqnarray}

Using $\mathbb{E}[P] = \theta/k$ and $\mathbb{E}[|X_1|^p ] = k^p
\mathbb{E}[| \Delta_k (F,Z_{k,1}) |^p ]$, and summing for
$k=1,\ldots,N$ we have proven \eqref{eqLyapunov}. By
condition~\eqref{eqinfvarcond}, if $p>2$, then this is
$o(V_N^{p/2})$ which by Lyapunov's theorem means that
$\frac{Y_N(F)- E_N}{\sqrt{V_n}}$ converges in law to a standard
Gaussian random variable as required.
\end{pf*}

\section*{Acknowledgment} A. Nikeghbali would like to thank G\'erard
Ben Arous for helpful discussions and for pointing out the use
of the improved bounds in the Feller coupling.

%

\printaddresses

\end{document}